\documentclass[11pt,twoside]{article} 
\usepackage[utf8]{inputenc}
\usepackage[paperwidth=165mm, paperheight=238mm,bottom=22mm,top=22mm,left=20mm,right=20mm]{geometry}
\usepackage{amsmath,amsthm,amssymb,amsfonts}
\usepackage{url}
\usepackage{color}
\usepackage{caption,subcaption}
\usepackage{textcomp}
\usepackage{rotating}
\usepackage{siunitx}
\usepackage{multirow}
\usepackage{float}
\usepackage{pdflscape}
\usepackage{afterpage}
\usepackage{cleveref}
\usepackage{soul}
\usepackage[linesnumbered, ruled]{algorithm2e}
\usepackage{pseudocode}
\usepackage[toc,page,titletoc,title]{appendix}

\newtheorem{theorem}{Theorem}[section]
\theoremstyle{plain}

\newtheorem{definition}[theorem]{Definition}
\newtheorem{lemma}[theorem]{Lemma}

\newtheorem{observation}[theorem]{Observation}
\newtheorem{question}[theorem]{Question}

\author{Gerold J\"ager \\ gerold.jager@umu.se 
	\and Klas Markstr\"om \\ klas.markstrom@umu.se 
	\and Denys Shcherbak \\ denys.shcherbak@umu.se
	\and Lars-Daniel \"Ohman \\ lars-daniel.ohman@umu.se 
	}

\date{}

\begin{document}

\title{Enumeration of Sets of Mutually Orthogonal Latin Rectangles}

\maketitle

\begin{abstract}
	We study sets of mutually orthogonal Latin rectangles (MOLR), and 
	a natural variation of the concept of self-orthogonal Latin squares 
	which is applicable on larger sets of mutually orthogonal Latin 
	squares and MOLR, namely that each Latin rectangle in a set of 
	MOLR is isotopic to each other rectangle in the set. We call such 
	a set of MOLR \emph{co-isotopic}.
	
	In the course of doing this, we perform a complete enumeration 
	of sets of $t$ mutually orthogonal $k\times n$ Latin 
	rectangles for $k\leq n \leq 7$, for all $t < n$ up to isotopism, 
	and up to paratopism. Additionally, for larger $n$ we
	enumerate co-isotopic sets of MOLR, as well as sets of MOLR
	where the autotopism group acts transitively on the  
	rectangles, and we call such sets of MOLR \emph{transitive}.
	
	We build the sets of MOLR row by row, and in this process we also 
	keep track of which of the MOLR are co-isotopic and/or transitive 
	in each step of the construction process. We use the prefix 
	\emph{stepwise} to refer to sets of MOLR with this property at each 
	step of their construction.
	
	Sets of MOLR are connected to other discrete objects, notably finite 
	geometries and certain regular hypergraphs. Here we observe that all 
	projective planes of order at most 9 except the Hughes plane can be 
	constructed from a stepwise transitive MOLR.
\end{abstract}

\section{Introduction}

Two Latin squares $L_A = (a_{ij})$ and $L_B = (b_{ij})$ of order $n$ 
are said to be \emph{orthogonal} if 
$| \{(a_{ij},b_{ij}): 1\leq i,j\leq n \}| = n^2$, that is, if when 
superimposing $L_A$ and $L_B$, we see each of the possible $n^2$ ordered
pairs of symbols exactly once. Pairs of orthogonal Latin squares show up in 
many different areas of combinatorics, both applied and pure.  
On the applied side, Latin squares are well known as the source of 
statistical designs for experiments. For more involved experiments 
one may instead switch to a design based on a pair of orthogonal 
Latin squares, or more generally one can use a set of pairwise 
orthogonal Latin squares, referred to as \emph{mutually orthogonal Latin 
squares} (MOLS). Here it is also possible to use mutually orthogonal 
$k \times n$ Latin rectangles (MOLR) when that suits the requirements 
of the statistical study. A less applied application of 
orthogonal Latin squares comes from the study of finite geometries. Here 
it is well known that a finite projective or affine plane of order $n$ 
exists if and only if a set of $(n-1)$ pairwise orthogonal $n \times n$ 
Latin squares exists. Motivated by connections like these, the first aim 
of this paper is to perform a complete enumeration of distinct, in the 
suitable sense, sets of $t$ mutually orthogonal $k \times n$ Latin 
rectangles, 
for as large values of the different parameters as possible.  Our second aim 
is to push this enumeration further for a few special classes of such 
sets of MOLR. Our computational results regarding MOLR here extend those 
by Harvey and Winterer in~\cite{HarWin}. Earlier work on MOLS is surveyed 
in the Handbook of Combinatorial Designs~\cite{Handbook}, and 
more recently,  Egan and Wanless~\cite{WanlessEgan} enumerated MOLS 
of order up to 9.

A particular type of orthogonal Latin squares which has seen use as a design 
is self-orthogonal Latin squares. A Latin square $L$ is said to be 
\emph{self-orthogonal} if $L$ and its transpose $L^T$ are orthogonal.  
The existence of self-orthogonal Latin squares (SOLS) of order $n$
for a number of values of $n$, not including $n=10$, was established
through many different constructions of infinite families of SOLS, 
e.g. in the seminal work by Mendelsohn~\cite{Mendelsohn71}.
In~\cite{Hedayat73}, Hedayat presented the first example of a self-orthogonal
Latin square of order $10$, and in~\cite{Brayton74}, 
Brayton et al. showed that $n\times n$ self-orthogonal Latin squares exist 
for each $n \neq 2,3,6$. Self-orthogonal Latin squares of orders up to and 
including $n = 9$ were completely enumerated by Burger et al. in 
~\cite{Burger10}, and by the same authors for order $10$ in~\cite{Burger10b}.

Being self-orthogonal is a special case of the concept of being 
\emph{conjugate orthogonal}, a concept introduced by Stein in~\cite{stein}. 
Here a conjugate is defined by a permutation $\sigma \in S_3$ which 
interchanges the roles of rows, columns, and symbols in the square. The 
transpose corresponds to the $\sigma$ which interchanges the roles of rows 
and columns. As Stein showed, for every $\sigma$ except the identity there 
are Latin squares which are orthogonal to their $\sigma$-conjugate, and 
Phelps~\cite{phelps} later investigated the possible orders of 
$\sigma$-conjugate Latin squares, settling the question with a handful of 
exceptions. These exceptions were later settled by Bennett, Wu and Zhu 
in~\cite{BWZ}.

Regarding conjugate orthogonality, one may sometimes find Latin squares which 
have several pairwise orthogonal conjugates, and it is possible to find 
examples where all six conjugates are pairwise orthogonal, as in 
Belyavskaya and Popovich~\cite{BP10} and Bennett~\cite{Ben}. 
This construction can of course not produce larger sets than six, 
and applying it to $k \times n$ Latin rectangles where $k<$ gives a 
maximum of four pairwise orthogonal rectangles.
However, if transposition is viewed as one of many possible `equivalence 
transformations' $\tau$, $\tau (L) = L^T$, then replacing $\tau$ with some 
other form of `equivalence' gives rise to a similar concept, where 
potentially a set of mutually orthogonal objects can be larger.

The main focus of the current paper is to enumerate sets of mutually 
orthogonal Latin rectangles. For small orders, we enumerate all MOLR
up to certain notions of equivalence, which we define in Section~\ref{sec_not}. 
For larger orders, where full enumeration of all MOLR is not 
feasible, we have proceeded using several natural subclasses of sets of MOLR. 
The first such class consists of those MOLR where each  
constituent Latin rectangle is isotopic to every other Latin rectangle in 
the set. We call such a set of MOLR \emph{co-isotopic}.
The next class is the \emph{transitive} sets of MOLR, where we call a set 
of MOLR transitive if the 
autotopism group of the set of MOLR acts transitively on the set of rectangles.  
So, every transitive set of MOLR is co-isotopic, but the reverse is not always 
true. Finally we consider \emph{stepwise transitive} and
\emph{stepwise co-isotopic} sets of MOLR, which are even more 
restricted classes, where we require that the set of MOLR can be constructed by 
adding one row at a time in such way that each set of MOLR along the way is 
transitive or co-isotopic, respectively. These conditions are very 
restrictive, but it turns out that all but one of the finite projective 
and affine planes of small orders have a corresponding stepwise 
transitive set of $(n-1)$ MOLS.

We will also discuss how, for the orders we have reached, our results lead 
to a complete enumeration of certain finite geometries, in particular, 
resolvable projective, affine and hyperbolic planes.

The work in the present paper builds on and extends results from a previous
paper of the current authors~\cite{triples}, where the focus was only on 
triples of MOLR. One distant goal is to approach the well-known long-standing 
open question of how large a maximum set of MOLS of order 10 is.

The paper is structured as follows. In Section~\ref{sec_not} we give basic 
notation and formal definitions regarding $t$-tuples of MOLR. In Section~\ref{sec_problem} 
we state the questions guiding our investigation, describe briefly the 
algorithm used to find all sets of MOLR and give some practical information 
regarding the computer calculations. In Section~\ref{sec_fingeom}, we give 
further background on finite geometries. In Section 5 we introduce a new reshaping transformation which maps a set of $t$ MOLR of size
$k\times n$ to a set of $(k-1)$ MOLR of size $(t+1)\times n$ and discuss 
some of the properties of this transformation.

In Section~\ref{sec_results} we present the data our computer search resulted 
in, which can be downloaded from~\cite{Web2}, together with further analysis 
and results. In particular, in Subsection~\ref{sec_totals} we give the total 
number of non-isotopic sets of MOLR for orders up to $n=7$. In Subsection~\ref{sec_larger} we discuss our enumeration 
of subclasses of sets of larger order, $n \geq 8$. Then in 
Subsection~\ref{sec_auto} we discuss the autotopism groups of the sets of MOLR. 

Finally, Section~\ref{sec_concl} concludes the main text with some open 
problems and observations. Here we also discuss the long-standing open 
problem of whether triples of mutually orthogonal Latin squares of order 10 
exist. In earlier works, Franklin~\cite{franklin84} constructed triples of 
MOLR, and Wanless~\cite{WANLESS2001} constructed 4-tuples of MOLR. We 
investigate the number of possible rows in stepwise transitive such examples.  
This is followed by a number
of appendices containing more detailed data on autotopism group sizes and 
all non-isotopic stepwise transitive $8$-MOLR with $n=9$. 

\section{Basic notation and definitions for $t$-MOLR}
\label{sec_not}

A \emph{Latin square} of order $n$ is an $n\times n$ matrix with cells 
filled by $n$ symbols, such that each row and each column contains each 
symbol exactly once. For $k \leq n$, a matrix with $k$ rows and $n$ 
columns whose cells are filled by $n$ symbols such that each row contains 
each symbol exactly once and each column contains each symbol at most 
once is called a $k\times n$ \emph{Latin rectangle}, and we shall refer
to $k\times n$ as its \emph{size}. In the following 
we use as symbol set $\{0,1,\dots,n-1\}$.

With \emph{mutually orthogonal Latin squares} defined as above, we can 
extend the orthogonality condition to Latin rectangles. We 
say that two $k\times n$ Latin rectangles $A=(a_{i,j})$ and $B=(b_{i,j})$ 
are \emph{orthogonal} if each ordered pair $(a_{i,j}, b _{i,j})$  
appears at most once. Also, a set of pairwise orthogonal Latin rectangles 
is called a set of \textit{mutually orthogonal Latin rectangles} (MOLR). 
A set of $t$ pairwise orthogonal Latin rectangles is called a $t$-MOLR 
for short.

We say that a $t$-MOLR is \emph{normalized} if it satisfies the 
following conditions:
\begin{itemize}
	\item[(S1)] (Ordering among columns) 
		The symbols in the first row of each rectangle appear in the order 
		$0, 1, \ldots, n-1$.
	\item[(S2)] (Ordering among rectangles) 
		The second row of the $i$:th rectangle is lexicographically larger 
		than the second row of the $(i+1)$:th rectangle. 
		In other words, if $a_1,a_2,\ldots,a_t$ are symbols in the 
		position $(2,1)$ in the $t$-MOLR, seen as an ordered $t$-tuple, 
		then it holds that $a_1>a_2> \ldots >a_t$.
	\item[(S3)] (Ordering among rows)
		The second row in the first rectangle is lexicographically larger 
		than the third row, the third row is larger than the fourth row, 
		and so on.
\end{itemize}
Note that (S2) relies on the rectangles having the same first row and 
being mutually orthogonal, so that all the entries in $(2,1)$ are in
fact distinct.

We use normalization to reduce the search space in
our computer runs, but note that two different normalized $t$-MOLR may
still be isotopic. For some further details, see our previous 
paper~\cite{triples}. As our computations proceed by adding consecutive 
rows to $t$-MOLR, we will have use for the following term: 
An \textit{extension} of size $k\times n$ is a $t$-MOLR which results 
from a $t$-MOLR of size $(k-1)\times n$ by adding one more row to each 
rectangle.

We will use several different notions of equivalence for $t$-MOLR. We 
recommend Section 2 of Egan and Wanless~\cite{WanlessEgan} for 
an in-depth discussion of the different symmetry and equivalence concepts 
for MOLS and we will follow that terminology. 
Let $(A_1,A_2,\ldots,A_t)$ be a $t$-MOLR of size $k\times n$, and let
$S_n$ denote the symmetric group on $n$ elements. The following group 
of isotopisms acts on the set of $t$-MOLR: $G_{n,k,t}=S_t\times S_k\times 
S_n \times [S_n \times S_n\times \ldots \times S_n]$, where the initial 
$S_t$ corresponds to permutations of the rectangles, the $S_k$ corresponds 
to permutations of the rows, the next $S_n$ corresponds to permutations 
of the columns, and each of the last $t$ copies of $S_n$, in square brackets,  corresponds to
permutations of the symbols in the single rectangles. 
Two $t$-MOLR $A$ and $A^\prime$ of size $k\times n$ are said to be \emph{isotopic} 
if there exists a $g\in G_{n,k,t}$ such that  
$g(A)=A^\prime$. The \emph{autotopism group} of a 
$t$-MOLR $A$ is defined as  $\mathrm{Aut}(A):=\{g\in G_{n,k,t}\, |\, g(A)=A\}$.

The autotopism group of a $t$-MOLR of size $k\times n$ is also a subgroup 
of a larger group, the \emph{autoparatopism} group, which can be described either using orthogonal 
arrays or, as we will do here,  via a hypergraph representation of the MOLR, which will also be used in one of our proofs.

\subsection{The standard hypergraph representation}

Let $A=(L^1,L^2,\ldots,L^t)$ be a list of $k\times n$ MOLR.
The \emph{standard labelled hypergraph representation $\mathcal{H}(A)$ of $A$} 
is constructed as follows.

The vertex set $V$ for $\mathcal{H}(A)$ consists of $t+2$ vertex classes, 
$V_1,V_2,\ldots V_{t+2}$. $V_1$ has one labelled vertex per row in the 
rectangles, and hence has size $k$. $V_2$ has one labelled vertex per column 
in the rectangles, and so has size $n$. For $i\geq 3$, $V_{i}$ has $n$ 
vertices  and for  each symbol from $L^{i-2}$ there is one vertex labelled 
with that symbol. 

The set of hyperedges of $\mathcal{H}(A)$ is as follows: For each $i\in V_1$ 
and $j\in V_2$ there is a hyperedge consisting of 
$(i,j,L^1_{i,j},L^2_{i,j},\ldots,L^t_{i,j})$.
Thus each hyperedge contains exactly one vertex from each vertex class, and each 
hyperedge corresponds to a position in the rectangles and the tuple of symbols 
used there by the rectangles.

The hypergraph $H(A)$ is $(t+2)$-uniform, and $(t+2)$-partite, with at most one 
class of vertices, $V_1$, smaller than $n$. This hypergraph is also linear, i.e. 
any pair of vertices belongs to at most one edge, since two edges intersecting 
in two vertices would either correspond to having more than one symbol in a 
given position $(i,j)$ in some rectangle, or a violation of orthogonality, or a 
repeated symbol in a row or column of some rectangle. 

Conversely, let $H$ be a linear $(t+2)$-uniform $(t+2)$-partite hypergraph 
with $t+1$ vertex classes $V_2,\ldots V_{t+2}$ of size $n$ and one vertex 
class $V_1$ of size $k\leq n$, and assume that the vertices in each class 
$V_i$ are labelled with a set of $|V_i|$ distinct symbols.  Given a choice 
of an ordered pair of vertex classes, which must include $V_1$ as the first 
class if $k<n$, we can construct a $t$-MOLR $A$ from $H$ by the same description 
as above, letting the first class be the row indices and the second class be 
the column indices. For that $t$-MOLR $A$ we have $H=\mathcal{H}(A)$ so we 
have an equivalence between $t$-MOLR and certain labelled hypergraphs. This 
is in fact just a different way of writing down the equivalence between 
$t$-MOLR and certain orthogonal arrays, with each hyperedge corresponding to 
a row in the orthogonal array.

\subsection{Isotopism and paratopism}

The isotopisms for a $k\times n$ $t$-MOLR $A$ correspond to the 
transformations which freely permute the labels of the vertices within each 
vertex class of $\mathcal{H}(A)$, and also permute the indices of the vertex 
classes $V_3,V_4,\ldots,V_{t+2}$.  

An isotopism is a special case of a more general type of transformation known 
as a \emph{paratopism}. Two $k\times n$ $t$-MOLR $A$ and $B$ are 
\emph{paratopic} if $\mathcal{H}(A)$ can be obtained from $\mathcal{H}(B)$ 
by permuting labels in the same way as for an isotopism, and freely permuting 
the indices of all vertex classes of equal size. The corresponding 
transformation of the MOLR is called a \emph{paratopism} on the set of MOLR.
Note that if $k<n$ then a paratopism between two MOLR cannot interchange the 
vertex class $V_1$ with any other class, since it is strictly smaller than the 
other classes, but for $k=n$ it may interchange $V_1$ with another class as well.  

For $k=n$, if the paratopisms are restricted to only interchanging $V_1$ with 
$V_2$, while permuting the other classes freely we get the \emph{trisotopisms}, 
and two $t$-MOLS related by a trisotopism are said to be \emph{trisotopic}.
Interchanging $V_1$ and $V_2$ corresponds to transposing the squares in the 
$t$-MOLS, so if two $t$-MOLS $A$ and $B$ are trisotopic they are either 
isotopic or the $A$ is isotopic to the $t$-MOLS obtained by transposing each 
square in $B$.

\section{Generation of $\mathbf{t}$-MOLR}
\label{sec_problem}

The basic method for our generation routine is quite simple. We start with 
the set of all $1 \times n$ $t$-MOLR  and find all possible extensions to  
$2 \times n$ $t$-MOLR, followed by an isotopy reduction where only one 
representative for each isotopy class is kept.  As part of the isotopy 
reduction, the autotopism group of each representative is also determined and 
its size stored. The extension step is then repeated until the desired number 
of rows is reached. The algorithms and methods used to generate the $t$-MOLR 
are rather straightforward extensions of those used in our previous 
paper~\cite{triples}. 
They were implemented in C++ and run in a parallelized version on the 
Kebnekaise and Abisko  supercomputers at High Performance Computing Centre 
North (HPC2N). The total run time for all the data in the paper was a few 
hundred core-years.

In order to safeguard the correctness of our computational results, several 
steps were taken. We wrote separate implementations in Mathematica of both 
the main algorithm and another much simpler method with which we performed 
an independent generation of the data for smaller sizes, in order to help 
verify the correctness of the C++ implementation. All the data has been compared 
with the known classifications of MOLS and MOLR in the literature, primarily 
McKay, Meynert and Myrvold~\cite{MMM} and Egan and Wanless~\cite{WanlessEgan}, 
and agrees with them.

The method described above was applied directly in order to first generate 
all $t$-MOLR of a given size when this was possible, and after each 
generation step we also classified the $t$-MOLR which belonged to one of
the following two classes.

\begin{definition}
\begin{itemize}
	\item[(a)]
	A $t$-MOLR $A=(A_1,A_2,\ldots,A_t)$ is \emph{co-isotopic} 
	if all pairs of rectangles $A_i,A_j\in A$ are isotopic.

	\item[(b)]
	A $t$-MOLR $A=(A_1,A_2,\ldots,A_t)$ is \emph{transitive} 
	if $\mathrm{Aut}(A)$ acts transitively on the set of rectangles in $A$. 
	That is, for all pairs $A_i,A_j\in A$ 
	there exists $\phi \in \mathrm{Aut}(A)$ which maps $A_i$ to $A_j$.
\end{itemize}
\end{definition}

When checking whether a $t$-MOLR is co-isotopic or transitive, we require 
the entire $t$-MOLR. This becomes a problem when generating $t$-MOLR by 
adding rows one by one. In order to find all co-isotopic $t$-MOLR we have 
to first generate all $t$-MOLR and then test them for co-isotopism and 
transitivity. We are therefore also interested in the two following 
recursively defined classes of $t$-MOLR which allow for more efficient 
generation.  

\begin{definition}\label{def_step}
\begin{itemize}
	\item[(a)]
	A $k \times n$ $t$-MOLR $A$ is \emph{stepwise co-isotopic} if $A$ is 
	co-isotopic and either $k=1$, or $k\geq2$  and $A$ is the extension of a 
	stepwise co-isotopic $(k-1) \times n$ $t$-MOLR. 

\item[(b)]
	A $k \times n$ $t$-MOLR $A$ is \emph{stepwise transitive} if $A$ is 
	transitive and either $k=1$, or $k\geq2$  and $A$ is the extension of a 
	stepwise transitive $(k-1) \times n$ $t$-MOLR. 
\end{itemize}
\end{definition}

As we shall demonstrate in Section~\ref{sec_larger}, the sets of 
$(n-1)$-MOLS corresponding to a classical projective plane over a finite 
field of order $n$ are, in fact, always stepwise transitive.

For both of the classes in Definition~\ref{def_step} we can generate the 
$k\times n$ $t$-MOLR by 
starting out with the corresponding class of $(k-1) \times n$ $t$-MOLR, 
finding all their non-isotopic extensions and then discarding those 
that are not co-isotopic or transitive, respectively. These restrictions 
typically lead to a far smaller set of $t$-MOLR to extend to the next value 
of $k$, and thanks to this we were able to perform complete enumeration of 
these classes for larger values of $n$ than in the general case.

In addition to producing the objects themselves, we also calculated the size 
of the autotopism group of each object.
With some exceptions due to size restrictions, all the data we generated is 
available for download at~\cite{Web2}. Further details about the organization 
of the data are given there.

We also include the size of the paratopism classes of MOLR in our tables. 
In order to obtain these classes, we used the package Nauty 
(see McKay and Piperno~\cite{MP14}) to reduce the set of isotopism classes 
to paratopism classes. In order to do this, the MOLR were encoded as graphs,
using the method described by Egan and Wanless in~\cite{WanlessEgan}.

\section{Finite geometries}
\label{sec_fingeom}

As mentioned in the introduction, sets of $n-1$ MOLS of size $n$ have a 
well known connection to finite projective and affine planes. Here we 
will recollect some facts from finite geometry and demonstrate how 
a general $t$-MOLR can be translated into the finite geometry setting. 

\subsection{Basic notions for finite geometries}
\label{basic_geom}

\begin{definition}
	A pair $\mathcal{P}=(V,L)$, where $V$ is a finite set of points and 
	$L$ is a set of subsets of $V$, which are called lines, is a 
	\emph{finite plane} if the following conditions are satisfied.
	\begin{enumerate}
		\item Each line has at least $2$ points.
		\item Any pair of points is contained in exactly one line.
		\item There exists a point $p$ and a line $\ell$, where $\ell$ 
			does not contain $p$.
		\item There exists a set of $4$ points such that no $3$ of 
			them lie on the same line.
	\end{enumerate}
	Additionally, the plane may satisfy one of the following parallelity
	properties:
	\begin{enumerate}
		\item[(P1)] Every pair of lines has non-empty intersection.
		\item[(P2)] Given a point $p$ and a line $\ell$, which does not contain 
			$p$, there exists exactly one line through $p$ which does not 
			intersect $\ell$.
		\item[(P3)] Given a point $p$ and a line $\ell$, which does not contain 
			$p$, there exist at least two lines through $p$ which do not 
			intersect $\ell$.
	\end{enumerate}
	A plane which satisfies $(P1)$ is called a \emph{finite projective plane}, 
	a plane which satisfies $(P2)$ is called a \emph{finite affine plane},  
	and a plane which satisfies $(P3)$ is called a \emph{finite hyperbolic plane}.
	
	If each line in a finite projective plane contains exactly $n+1$ points, the
	finite projective plane is said to be of \emph{order} $n$, and correspondingly,
	if each line in a finite affine plane contains exactly $n$ points, the finite 
	affine plane is said to be of order $n$.
\end{definition}

A collection $P$ of non-intersecting lines from $\mathcal{P}$ which form 
a partition of $V$ is called a \emph{parallel class}.  A partition of the 
lines of $\mathcal{P}$ into parallel classes is called a \emph{resolution} 
of $\mathcal{P}$, and a geometry which has at least one resolution is 
called a \emph{resolvable} geometry.   

A set of $n-1$ MOLS of order $n$ is equivalent to a finite projective plane 
of order $n$, as pointed out by Bose~\cite{constr}. The existence of 
a set of $n-1$ MOLS, for $n=p^r$, for a prime $p$, was demonstrated in a 
somewhat forgotten paper from \num{1896} by Moore~\cite{Moo96}. 
See Ehrhardt~\cite{EHR} for a historical discussion of Moore's paper and 
the wider context of 19th-century design theory.

To relate our results on Latin rectangles and squares to the context of finite 
geometries, we will first recall the explicit correspondence between finite 
projective geometries on the one hand, and complete sets of mutually orthogonal 
Latin squares on the other hand. We will begin with noting that a finite 
projective plane of order $n$ has $n^2+n+1$ points and $n^2+n+1$ lines, 
each containing $n+1$ points, and that each point in a finite projective plane 
of order $n$ belongs to $n+1$ lines.

To get the correspondence between a finite projective plane of order $n$ and 
a set of $n-1$ MOLS of order $n$, we first select a line of the projective 
plane, which we call the \emph{line at infinity}, $L_\infty$. Next, we select
two points $x_\infty$ and $y_\infty$ on $L_\infty$, and label the remaining 
$n-1$ points on $L_\infty$ by $\ell_1, \ell_2, \ldots, \ell_{n-1}$. The $n$ 
remaining lines containing $x_\infty$ are labelled by $X_1,X_2,\ldots,X_n$, 
and the $n$ remaining lines containing $y_\infty$ are labelled
$Y_1, Y_2, \ldots, Y_n$.

The points $\ell_1, \ell_2, \ldots, \ell_{n-1}$ will correspond to the $n-1$ 
Latin squares $L_1, L_2, \ldots L_{n-1}$ in the set of MOLS, and the point 
$p_{i,j}$ at the intersection between $X_i$ and $Y_j$ will correspond to the 
cell $(i,j)$ in the Latin squares. The symbol that goes in cell $(i,j)$ in 
$L_k$ is given by fixing a labelling using symbols $s_1, s_2, \ldots,s_n$ of 
the $n$ lines (excluding $L_\infty$) that contain $\ell_k$, and checking 
which of the symbols was assigned to the unique line passing through $p_{i,j}$.

Note that the ample room for choices of labelings gives different 
sets of MOLS, and that the construction can be reversed, that is, any 
$(n-1)$-MOLS gives a finite projective plane.


When $\mathcal{P}_n$ is the classical Galois plane over the finite field 
$GF(n)$, sometimes denoted by $PG(2,n)$, one can give a simple explicit form 
for a set of MOLS derived from that plane in the following way. 
Let $x\in GF(n)$ be a generator for the multiplicative group of $GF(n)$ 
and set $a_i=x^{i-1}$ for $i=1,\ldots, n$.
Now, for $1\leq k \leq n-1$ define a Latin square $L_k$ by setting position 
$(i,j)$ equal to $a_i+ a_k\times a_j$. As shown, e.g., in D\'enes and Keedwell's 
Latin squares book~\cite{DK15}, this defines an $(n-1)$-MOLS, which 
in turn also defines the classical Galois plane over $GF(n)$. 


When $n$ is a prime power, finite projective planes can always be constructed 
using the Galois field $GF(n)$, but there are projective planes that do not 
arise in this way. For example, as early as 1907, 
Veblen and MacLagan-Wedderburn~\cite{Veblen1907}  
constructed $3$ projective planes of order~$9$, not isomorphic to the standard 
projective plane of order $9$ arising from $GF(9)$. 
Considerably later, Lam, Kolesova and Thiel~\cite{Lam91} showed by an exhaustive computer 
search that these $4$ projective planes of order $9$ are, in fact, the only ones.

For some other $n$, notably $n=6$ and $n=14$, the Bruck-Ryser 
theorem~\cite{bruckryser} excludes the possibility of a projective plane.
For $n=6$, non-existence was already known, since Tarry~\cite{tarry} had proven 
that Euler's 36 officers problem (which asked for a set of just two MOLS of order 
$6$) had no solution, leaving $n=10$ the smallest open case. A delightful account 
of the search for a projective plane of order~$10$ can be found in 
Lam~\cite{Lam}, and the non-existence was settled by Lam, Thiel and 
Swiercz in~\cite{LTS}. Combining this computational 
non-existence result with a result of Shrikhande~\cite{Shrikhande61}, one gets 
that there does not exist a set of $7$ or more MOLS of order $10$. Since there 
are examples of pairs of orthogonal Latin squares of order $10$, the maximum 
order of a $t$-MOLS for $n=10$ lies in the interval $2 \leq t \leq 6$.

Given a projective plane $\mathcal{P}$ of order $n$ we can construct an 
affine plane of the same order by deleting a line and all the points on it 
from $\mathcal{P}$. It is also well known that any affine plane can be obtained 
from a projective plane in this way. So, the existence of an affine  or a 
projective plane of a given order are equivalent, and in turn equivalent to the 
existence of a set of $n-1$ MOLS of order $n$. Finite hyperbolic geometries 
have not been studied in as great detail as the projective or affine ones. 
The first axiomatization 
for finite hyperbolic planes was given by Graves in~\cite{Grav} and a 
few early constructions and structural theorems were given 
by Di Paola, Henderson, and Ostrom in~\cite{DiPa,Hen62,Ost}, respectively. 
Sandler~\cite{San} also noted that one will obtain a 
finite hyperbolic plane from a finite projective plane by deleting three 
lines that
do not intersect in a single point, together with all points on these lines, 
or equivalently by deleting two non-parallel lines from a finite affine plane.

Another well studied class of finite geometries are \emph{nets}. They were 
introduced by Bruck~\cite{Bru1,Bru2} who also showed that they, analogous 
to the situations for projective and affine planes, are equivalent to 
$t$-MOLS. A net $\mathcal{P}(A)$ can be constructed in the following way. 
Let $A=(A^1,A^2,\ldots,A^t)$ be a set of mutually orthogonal Latin squares
of order $n$, and let the set of pairs 
$V=\{(i,j) \mid 1\leq i \leq n,1\leq j\leq n\}$ be the 
point set. For each square $A^s$ and symbol $r$ we let the set of points 
such that $A^s_{i,j}=r$ be a line in $L$. We also add one line to $L$ for 
each row of the squares and one for each column, and set $\mathcal{P}(A)=(V,L)$.
We now note that each line in $\mathcal{P}(A)$ has $n$ points, each point lies 
in $t+2$ lines, two lines defined by different squares intersect in exactly 
one point (since the squares are orthogonal), the set of lines defined by a 
single square are pairwise disjoint, and any pair of points lies in at most  
one line (again by orthogonality). 
Additionally, the set of lines defined by a single square form a partition of 
$V$, i.e., a parallel class, and the parallel classes defined by the different 
squares form a partition of $L$, so $\mathcal{P}(A)$ is a resolvable geometry.     
The net $\mathcal{P}(A)$ is of order $n$ and degree $t$. Note that 
if $t=n-1$, then the net is an affine plane, and in general, nets are a 
particular class of \emph{partial geometries} as defined by Bose~\cite{Bose63}.     

\subsection{Finite geometries from Latin rectangles}
\label{rec_geom}

The classical finite geometric constructions which we have surveyed here 
are all connected to Latin squares. However, many of these constructions 
can fruitfully be extended to Latin rectangles and MOLR as well and we 
will briefly discuss how this can be done 

If we apply the same construction as for nets to a $k \times n$ $t$-MOLR 
we get a weaker geometry which we call a \emph{partial net}. When $k<n$ 
we lose the property that every pair of lines from different parallel 
classes intersect. Of particular interest here are $t$-MOLR which are 
maximal with respect to either $t$ or $k$, since these give rise to 
geometries which cannot be embedded in a larger partial net with the same 
value of $n$.

Our more symmetric classes of $t$-MOLR correspond to geometries with certain 
symmetries. A transitive $t$-MOLR gives rise to a geometry with the property 
that the automorphism group of the geometry acts transitively on the set of 
parallel classes in a specific resolution of the geometry. A stepwise 
transitive set in turn gives us a geometry with the property that the 
geometry can be built by stepwise adding new points in a way which preserves 
the transitivity of the resolution.

\section{A reshaping transformation}
\label{sec_reshape}

For the specific case of $(n-1)$-MOLS, several authors (see Egan and 
Wanless~\cite{WanlessEgan} for further details) have studied a notion of 
equivalence where two such MOLS are equivalent if they define the same 
projective plane. Using the relation between MOLS and projective planes 
described in Section~\ref{sec_fingeom}, this means that two such 
equivalent MOLS are constructed by selecting different lines at infinity 
from the projective plane. For $t$-MOLR there is in general no extension 
to a projective plane, but we will define a similar transformation $T$ 
which maps a $k\times n$ $t$-MOLR $A$ to a $(t+1)\times n$ $(k-1)$-MOLR 
$T(A)$.

In preparation of defining this transformation, we will first define a 
hypergraph $G(A) = (W,E)$ with labelled vertices and hyperedges, from a
$k\times n$ $t$-MOLR $A=\{L^1,L^2,\ldots,L^t \}$ where each 
rectangle uses the symbols $1,2,\ldots,n$. The hypergraph $G(A)$ is similar, 
but not identical, to the standard hypergraph representation.

Let the vertex set $W$ be the set of ordered pairs $(i,j)$, where 
$1\leq i \leq k$ and $1\leq j \leq n$ together with the singletons 
$0\leq q \leq  t$. A vertex $(i,j)$ is given the label $i$ and a singleton 
vertex $q$ is labelled by the corresponding symbol $q$.

Next we will define the set $E$ of labelled hyperedges. For each $q\geq 1$ 
and symbol $s$ we add a hyperedge consisting of the singleton 
$q$ and the pairs $(i,j)$ such that $L^q_{ij}=s$, and label this edge by $s$.
We think of the singletons $q$ as an indexing of the rectangles, and the
edges containing the singleton $q$ keep track of symbol positions in $L^q$.
Permuting symbols in the Latin rectangles $L^1, L^2, \ldots L^t$ thus only
affects the edge labels of $G$.
For each row index $1\leq i \leq k$ we add a labelled hyperedge which 
consists of all pairs $(i,j)$ with $j=1,2,\ldots, n$ and the singleton $q=0$, 
labelled by $i$. 

Note that for $k=n$, $t=n-1$ the labelled hypergraph $G(A)$, 
with labels removed and interpreting edges as lines, is equivalent to the 
projective plane defined by $A$ with the point in the line at infinity 
where the column-lines intersect removed, and the column-lines deleted.

\begin{lemma}\label{Glin}
	The hypergraph $G(A)=(W,E)$ is linear, $(k+1)$-uniform, and 
	$(k+1)$-partite, with $k$ vertex classes $W_1, W_2,\ldots, W_k$ 
	of size $n$ corresponding to the rows of $A$, and one vertex class 
	$W_0$ of size $t+1$ corresponding to the added singletons.
\end{lemma}

\begin{proof}
	By construction, $W_i$, $1\leq i\leq k$, has size $n$,
	and $W_{0}$ has size $t+1$. Each hyperedge contains 
	exactly one vertex from each vertex class, so the hypergraph is 
	$(k+1)$-uniform and $(k+1)$-partite.
	
	In $G(A)$, two edges intersecting in 
	a singleton $q$ have no further common vertices since they either 
	correspond to positions of different symbols in the same rectangle, or 
	to distinct columns of a rectangle. Edges containing distinct singletons 
	cannot intersect in two other vertices, since that would either mean 
	that the same symbol appears twice in a column for the same rectangle, 
	or that a pair of symbols appears twice in a pair of rectangles, hence 
	violating orthogonality. So the hypergraph is linear.
\end{proof}

Lemma~\ref{Glin} says that $G(A)$ satisfies the properties required for 
being the standard hypergraph representation for some MOLR. Hence we can 
interpret $G(A)$, with a specified $(k+1)$-partition of its vertex set, 
as the standard hypergraph representation $\mathcal{H}(T(A))$ for some 
$(k-1)$-tuple $T(A)$ of $(t+1)\times n$ MOLR, by interpreting the vertex 
class $W_{0}$ as corresponding to $V_1$, the row indices for $T(A)$, and 
taking one of the other $W_i$, which we take to be $W_1$, as corresponding 
to $V_2$, the column indices in $T(A)$. The remaining $W_2, W_3, \ldots, W_k$
are interpreted as the indexing of the rectangles, $V_i = W_{i-1}$ for 
$i = 3, \ldots, k+1$.

This implicitly defines a transformation $T$ from a $k\times n$ 
$t$-MOLR $A$ to a $(t+1)\times n$ $(k-1)$-MOLR $T(A)$, so if one of the MOLR 
exists, then so does
the other. Note that $T(T(A))$ gives a MOLR with the same parameters as $A$. 
However, the transformation does not induce an isomorphism of the paratopism 
group for the set of $k\times n$ $t$-MOLR to that for $(t+1)\times n$ 
$(k-1)$-MOLR. For most parameters these two paratopism groups do not have 
the same size, snd so cannot be isomorphic. This in turn means that we do 
not necessarily have the same number of paratopism classes for $k\times n$ 
$t$-MOLR and $(t+1)\times n$ $(k-1)$-MOLR. For small $n$, with $t<n-1$, it 
turns out from our computational results that the number of paratopism 
classes coincide, see Tables~\ref{MOLR4p}, \ref{MOLR5p}, \ref{MOLR6p} and 
\ref{MOLR7p}, but we expect this to fail for $n$ slightly larger than those 
we consider here. A full investigation of how the transformation $T$ 
interacts with paratopism is beyond the scope of the current article, but 
we note some of the properties of $T$ in the following theorem.

\begin{theorem}
	The mapping $T$ maps the set of $k\times n$ $t$-MOLR 
	to the set of $(t+1)\times n$ $(k-1)$-MOLR.
	
	Any MOLR $A$ with the identity permutation as the first row in each 
	rectangle is mapped by $T$ to another MOLR with the identity permutation 
	as the first row in each rectangle.
	
	Two MOLR $A$ and $B$ have $T(A)=T(B)$ if and only if $B$ can be 
	obtained from $A$ by only permuting the symbols of each rectangle. 
\end{theorem}

\begin{proof}
	It follows directly from the definition that $T$ maps the set of 
	$k\times n$ $t$-MOLR to the set of $(t+1)\times n$ $(k-1)$-MOLR.
	Since each edge in $G(A)$ incident to the vertex given by the singleton 
	0 consists of the pairs $(i,j)$ with a fixed second coordinate, the MOLR 
	$T(A)$ has the identity permutation as the first row of each rectangle 
	if $A$ does.
	
	Suppose $B$ can be obtained from $A$ by only permuting symbols in the
	rectangles. Then $G(A)$ and $G(B)$ with removed edge labels will leave the
	same unlabelled hypergraph, and thus represent the same MOLR, up to 
	permutation of symbols in each rectangle. As we have a fixed choice of 
	index classes $W_0 = V_1$,  $W_1 = V_2$, this in turn means that $T(A)=T(B)$.
	
	Conversely, suppose $T(A) = T(B)$. Then the standard hypergraph 
	representations of $T(A)$ and $T(B)$ of course coincide,
	$\mathcal{H}(T(A))=\mathcal{H}(T(B)) = \mathcal{H}$.
	Removing the edge labels from $G(A)$ and $G(B)$ leaves the same 
	unlabelled hypergraph $\mathcal{H}$ if $G(A)$ and $G(B)$ differed only 
	in the order of the edge labels, corresponding to permutations of the 
	symbols in some of the rectangles of $A$ and $B$. 
\end{proof}

%

\section{Computational results and analysis}
\label{sec_results}

We now turn to the results and analysis of our computational work. 

\subsection{The number of $\mathbf{t}$-MOLR}
\label{sec_totals}

In Tables~\ref{MOLR4}--\ref{MOLR7}, 
we present data on the number of isotopism classes and paratopism classes of 
$k \times n$ $t$-MOLR. 

\begin{table}[H]
	\begin{center}
		\begin{tabular}{|r||r|r|r|}
			\hline 
			&	\multicolumn{1}{c|}{2$\times$4}   & \multicolumn{1}{c|}{3$\times$4}   & \multicolumn{1}{c|}{4$\times$4}   \\ 
			\hline \hline
			$t=1$	&	2	&	2	&	2	 \\ \hline 
			$t=2$	&	3	&	2	&	1	 \\ \hline 
			$t=3$	&	2	&	1	&	1	 \\ \hline 	
	\end{tabular}\end{center}
	\caption{The number of isotopism classes of $t$-MOLR for $n=4$.}
	\label{MOLR4}
\end{table}

\begin{table}[H]
	\begin{center}
		\begin{tabular}{|r||r|r|r|}
			\hline 
			&	\multicolumn{1}{c|}{2$\times$4}   & \multicolumn{1}{c|}{3$\times$4}   & \multicolumn{1}{c|}{4$\times$4}   \\ 
			\hline \hline
			$t=1$	&	2	&	2	&	2	 \\ \hline 
			$t=2$	&	2	&	2	&	1	 \\ \hline 
			$t=3$	&	2	&	1	&	1	 \\ \hline 	
	\end{tabular}\end{center}
	\caption{The number of paratopism classes of $t$-MOLR for $n=4$.}
	\label{MOLR4p}
\end{table}

\begin{table}[H]
	\begin{center}
		
		\begin{tabular}{|r||r|r|r|r|}
			\hline 
			&	\multicolumn{1}{c|}{2$\times$5}   & \multicolumn{1}{c|}{3$\times$5}   & \multicolumn{1}{c|}{4$\times$5}   & \multicolumn{1}{c|}{5$\times$5}   \\ 
			\hline \hline 
			$t=1$	&	2	&	3	&	3	&	2	 \\ \hline 
			$t=2$	&	5	&	14	&	2	&	2	 \\ \hline 
			$t=3$	&	4	&	1	&	1	&	1	 \\ \hline 
			$t=4$	&	3	&	1	&	1	&	1	 \\ \hline 
	\end{tabular}\end{center}
	\caption{The number of isotopism classes of $t$-MOLR for $n=5$.}
	\label{MOLR5}
\end{table}

\begin{table}[H]
	\begin{center}
		
		\begin{tabular}{|r||r|r|r|r|}
			\hline 
			&	\multicolumn{1}{c|}{2$\times$5}   & \multicolumn{1}{c|}{3$\times$5}   & \multicolumn{1}{c|}{4$\times$5}   & \multicolumn{1}{c|}{5$\times$5}   \\ 
			\hline \hline 
			$t=1$	&	2	&	3	&	3	&	2	 \\ \hline 			
			$t=2$	&	3	&	9	&	1	&	1	 \\ \hline 
			$t=3$	&	3	&	1	&	1	&	1	 \\ \hline 
			$t=4$	&	2	&	1	&	1	&	1	 \\ \hline 
	\end{tabular}\end{center}
	\caption{The number of paratopism classes of $t$-MOLR for $n=5$.}
	\label{MOLR5p}
\end{table}

\begin{table}[H]
	\begin{center}
		
		\begin{tabular}{|r||r|r|r|r|r|}
			\hline 
			&	\multicolumn{1}{c|}{2$\times$6}   & \multicolumn{1}{c|}{3$\times$6}   & \multicolumn{1}{c|}{4$\times$6}   & \multicolumn{1}{c|}{5$\times$6}   & \multicolumn{1}{c|}{6$\times$6}   \\ 
			\hline \hline
			$t=1$	&	4	&	\num{16}	&	\num{56}	&	40	&22 \\ \hline 			
			$t=2$	&	28	&	\num{1526}	&	\num{2036}	&	85	&0 \\ \hline 
			$t=3$	&	103	&	\num{2572}	&	513	&	7 &0	 \\ \hline 
			$t=4$	&	92	&	118	&	12	&	8 &0	 \\ \hline 
			$t=5$	&	33	&	0	&0 &0 &0 \\ \hline 
			
	\end{tabular}\end{center}
	\caption{The number of isotopism classes of $t$-MOLR for $n=6$.}
	\label{MOLR6}
\end{table}

\begin{table}[H]
	\begin{center}
		
		\begin{tabular}{|r||r|r|r|r|r|}
			\hline 
			&	\multicolumn{1}{c|}{2$\times$6}   & \multicolumn{1}{c|}{3$\times$6}   & \multicolumn{1}{c|}{4$\times$6}   & \multicolumn{1}{c|}{5$\times$6}   & \multicolumn{1}{c|}{6$\times$6}   \\ 
			\hline \hline
			$t=1$	&	4	&	\num{14}	&	\num{44}	&	33	&12 \\ \hline 						
			$t=2$	&	14	&	\num{575}	&	\num{745}	&	44	&0 \\ \hline 
			$t=3$	&	44	&	\num{745}	&	179	        &	5 &0	 \\ \hline 
			$t=4$	&	33	&	44	         &	5       	&	5 &0	 \\ \hline 
			$t=5$	&	17	&	0	&0 &0 &0 \\ \hline 
			
	\end{tabular}\end{center}
	\caption{The number of paratopism classes of $t$-MOLR for $n=6$.}
	\label{MOLR6p}
\end{table}

\begin{table}[H]
	\begin{center}
		\setlength\tabcolsep{3pt}
		\begin{tabular}{|r||r|r|r|r|r|r|}
			\hline 
			&	\multicolumn{1}{c|}{2$\times$7}   & \multicolumn{1}{c|}{3$\times$7}   & \multicolumn{1}{c|}{4$\times$7}   & \multicolumn{1}{c|}{5$\times$7}   & \multicolumn{1}{c|}{6$\times$7}   & \multicolumn{1}{c|}{7$\times$7}   \\ 
			\hline \hline
			$t=1$	&	4	&	\num{56}	&	\num{1398}	&	\num{6941}	&	\num{3479}	&	564	 \\ \hline 			
			$t=2$	&	100	&	\num{514162}	&	\num{49415812}	&	\num{21290125}	&	\num{11582}	&	20	 \\ \hline 
			$t=3$	&	 2858 &	\num{65883453}	&	\num{323112477} 	&	\num{55545}	&	16	&	4	 \\ \hline 
			$t=4$	&	\num{17609}	&	\num{35469948}	&	\num{68659}	&	204	&	7	&	3	 \\ \hline 
			$t=5$	&	\num{10626}	&	\num{22982}	&	19	&	5	&	5	&	1	 \\ \hline 
			$t=6$	&	\num{1895}	&	23	&	2	&	1	&	1	&	1	 \\ \hline 
			
	\end{tabular}\end{center}
	\caption{The number of isotopism classes of $t$-MOLR for $n=7$.}
	\label{MOLR7}
\end{table}

\begin{table}[H]
	\begin{center}
		\setlength\tabcolsep{3pt}
		\begin{tabular}{|r||r|r|r|r|r|r|}
			\hline 
			&	\multicolumn{1}{c|}{2$\times$7}   & \multicolumn{1}{c|}{3$\times$7}   & \multicolumn{1}{c|}{4$\times$7}   & \multicolumn{1}{c|}{5$\times$7}   & \multicolumn{1}{c|}{6$\times$7}   & \multicolumn{1}{c|}{7$\times$7}   \\ 
			\hline \hline
			$t=1$	&	4	&	\num{45}	&	\num{808}	&	\num{3712}	&	\num{1895}	&	147	 \\ \hline 						
			$t=2$	&	45	&	\num{172622}	&    \num{16481351}	&	\num{7103198}	&	\num{4013}	&	7	 \\ \hline 
			$t=3$	&	 808 &	\num{16481351} & \num{80797488 } 	&	\num{14121}	&	12	&	1	 \\ \hline 
			$t=4$	&	\num{3712}       	    &	\num{7103198}	&	\num{14121}	&	82	&	4	&	1	 \\ \hline 
			$t=5$	&	\num{1895}	            & 	\num{4013}	&	12	&	4	&	4	&	1	 \\ \hline 
			$t=6$	&	\num{324}	&	11	   &	2               	&	1	&	1	&	1	 \\ \hline 
			
	\end{tabular}\end{center}
	\caption{The number of paratopism classes of $t$-MOLR for $n=7$.}
	\label{MOLR7p}
\end{table}


In the data for $n \leq 7$ some patterns can be observed, somewhat 
interrupted by the exceptional behavior for $n=6$. If we consider 
fixed values of $t$ and $n$ and increase $k$ we always see a unimodal 
sequence, and the peak of the sequence appears at a lower value of 
$k$ when $t$ is increased. The patterns conform well with the number 
of constraints on the symbols, as a function of $t$ and $k$. If we 
instead keep $n$ and $k$ fixed and increase $t$ we see a similar pattern, 
though here there is an exception for $n=6$, $k=5$, where there is a 
local minimum at $t=3$. These observations motivate the following 
questions.

\begin{question}
	Is the number of $t$-MOLR for fixed $n$ and $t$ a unimodal sequence in $k$?  
\end{question}

\begin{question}
	For $n\geq7$, is the number of $t$-MOLR for fixed $n$ and $k$ a unimodal 
	sequence in $t$?  
\end{question}

Additionally we see in the tables for paratopism classes that there is a 
symmetry between the entries below and above the diagonal, with exceptions 
for the cases with $k=n$ or $t=n-1$. For the case of $6 \times 6$ Latin 
squares (that is, $t=1$, $k=n=6$), this exception corresponds to a change from 
paratopism to trisotopism as the equivalence relation, 17 is the number of 
trisotopism classes and 12 the number of paratopism classes (see Egan and 
Wanless~\cite{WanlessEgan}). For $n=7$ we have three entries for $k=n$ which 
do not match the corresponding ones for $t=n-1$. For these small values 
of $n$ this partial symmetry is explained by the reshaping transformation 
described in Section~\ref{sec_reshape}, but as discussed there we expect 
the symmetry to fail in general for large enough $n$.

We have also classified the small sets of MOLR according to some further 
properties. In Tables~\ref{4x4.}--\ref{6x6.}, for $n=4,5,6$, we give the 
number of $t$-MOLR that are A) co-isotopic, B) transitive, C) stepwise 
co-isotopic and D) stepwise transitive. In a sense, these four classes 
are gradually more regular, and the data in the tables gives the total 
numbers from each such class in the form A, B, C, D in each cell. In 
Table~\ref{7x7}, the data is presented in the form A, B, and in 
Table~\ref{7x7STEP}, the data is in the form C, D.
Comparing the data in Tables~\ref{MOLR4}-\ref{MOLR7} with the data in 
Tables~\ref{4x4.}-\ref{7x7STEP}, it is clear that when $k$ or $t$ is small 
compared to $n$, most $t$-MOLR have none of these stronger regularity 
properties, but whenever a $t$-MOLR exists, we also have a co-isotopic 
$t$-MOLR with the same parameters. 

Whenever there exists a $t$-MOLR, we also find transitive $t$-MOLR for most 
parameters. The exception is $n=7$, where there exist $t$-MOLS ($k=7$, that is) 
with $t=4,5$ but no corresponding transitive $t$-MOLR, demonstrating that the 
autotopism group for the $6$-MOLR does not have orbits of length 4 and 5. 
As a further example of observations from the data, for $n=7, k=4$ 
(see Table~\ref{7x7}) there exist co-isotopic $5$-MOLR, but no transitive 
$5$-MOLR, and \emph{a fortiori}, no stepwise transitive $5$-MOLR.

\begin{table}[H]
	\begin{center}
		\begin{tabular}{|c||c|c|c|}
			\hline 
			&	{2$\times$4}   & {3$\times$4}   & {4$\times$4}   \\ 
			\hline \hline
			$t=2$	&	2, 2, 2, 2	&	2, 2, 1, 1	&	1, 1, 1, 1	 \\ \hline 
			$t=3$	&	1, 1, 1, 1	&	1, 1, 1, 1	&	1, 1, 1, 1	 \\ \hline 
	\end{tabular}\end{center}
	\caption{The number of non-isotopic $t$-MOLR for $n=4$ sorted by increasing regularity.}
	\label{4x4.}
\end{table}

\begin{table}[H]
	\begin{center}
		
		\begin{tabular}{|c||c|c|c|c|}
			\hline 
			&	{2$\times$5}   & {3$\times$5}   & {4$\times$5}   & {5$\times$5}   \\
			\hline \hline
			$t=2$	&	4, 3, 4, 3	&	11, 9, 7, 6	&	2, 1, 2, 1	&	2, 1, 2, 1	 \\ \hline 
			$t=3$	&	3, 2, 3, 2	&	1, 0, 0, 0	&	1, 0, 0, 0	&	1, 0, 0, 0	 \\ \hline 
			$t=4$	&	2, 2, 2, 2	&	1, 1, 1, 1	&	1, 1, 1, 1	&	1, 1, 1, 1	 \\ \hline 
			
	\end{tabular}\end{center}
	\caption{The number of non-isotopic $t$-MOLR for $n=5$ sorted by increasing regularity.}
	\label{5x5.}
\end{table}

\begin{table}[H]
	\begin{center}
		\setlength\tabcolsep{2pt}
		\begin{tabular}{|c||c|c|c|c|}
			\hline 
			&	{2$\times$6}   & {3$\times$6}   & {4$\times$6}   & {5$\times$6} \\ 
			\hline \hline
			$t=2$	&	12, 11, 12, 11	&	280, 170, 158, 103	&	229, 160, 66, 50	&	43, 36, 13, 12	 \\ \hline 
			$t=3$	&	16, 6, 16, 6	&	115, 29, 32, 4	&	62, 39, 4, 1	&	4, 3, 0, 0	 \\ \hline 
			$t=4$	&	9, 8, 9, 8	&	19, 17, 15, 15	&	4, 3, 0, 0	&	4, 4, 0, 0	 \\ \hline 
			$t=5$	&	2, 2, 2, 2	&	0, 0, 0, 0 &	0, 0, 0, 0 &	0, 0, 0, 0	 \\ \hline 
			
	\end{tabular}\end{center}
	\caption{The number of non-isotopic $t$-MOLR for $n=6$ sorted by increasing regularity.}
	\label{6x6.}
\end{table}

\begin{table}[H]
	\setlength\tabcolsep{2pt}
	\begin{center}
		\begin{tabular}{|c||c|c|c|c|c|c|}
			\hline 
			&	{2$\times$7}   & {3$\times$7}   & {4$\times$7}   & {5$\times$7}   & {6$\times$7}   & {7$\times$7}   \\ 
			\hline \hline
			$t=2$	&	42, 29	&	\num{14464}, \num{3549}	&	\num{65156}, \num{27299}	&	\num{22432}, \num{18836}	&	409, 392	&	9, 6	 \\ \hline 
			$t=3$	&	318, 15	&	\num{49370}, 647	&	\num{2985}, \num{1578}	&	111, 36	&	11, 6	&	4, 1	 \\ \hline 
			$t=4$	&	691, 21	&	\num{1622}, 110	&	84, 67	&	67, 53	&	7, 3	&	3, 0	 \\ \hline 
			$t=5$	&	176, 6	&	49, 42	&	2, 0	&	4, 2	&	5, 3	&	1, 0	 \\ \hline 
			$t=6$	&	26, 5	&	7, 7	&	2, 2	&	1, 1	&	1, 1	&	1, 1	 \\ \hline 
			
		\end{tabular}
	\end{center}
	\caption{The number of non-isotopic co-isotopic and transitive $t$-MOLR for $n=~7$.}
	\label{7x7}
\end{table}

\begin{table}[H]
	\setlength\tabcolsep{2pt}
	\begin{center}
		\begin{tabular}{|c||c|c|c|c|c|c|}
			\hline 
			&	{2$\times$7}   & {3$\times$7}   & {4$\times$7}   & {5$\times$7}   & {6$\times$7}   & {7$\times$7}   \\ 
			\hline \hline
			$t=2$	&	42, 29	&	\num{7423}, \num{2175}	&	\num{14960}, \num{10029}	&	\num{4163}, \num{3923}	&	91, 84	&	6, 4	 \\ \hline 
			$t=3$	&	318, 15	&	\num{13975}, 185	&	283, 160	&	8, 5	&	4, 1	&	4, 1	 \\ \hline 
			$t=4$	&	691, 21	&	585, 48	&	12, 1	&	3, 0	&	3, 0	&	3, 0	 \\ \hline 
			$t=5$	&	176, 6	&	48, 42	&	2, 0	&	1, 0	&	1, 0	&	1, 0	 \\ \hline 
			$t=6$	&	26, 5	&	6, 4	&	2, 2	&	1, 1	&	1, 1	&	1, 1	 \\ \hline  
			
		\end{tabular}
	\end{center}
	\caption{The number of non-isotopic stepwise co-isotopic and stepwise transitive $t$-MOLR for $n=7$.}
	\label{7x7STEP}
\end{table}

The numbers of stepwise co-isotopic and stepwise transitive $t$-MOLR are by 
definition smaller than (or equal to) the numbers of co-isotopic and transitive 
$t$-MOLR respectively, but we again find stepwise co-isotopic examples for most 
parameter values. For $n=7$, $k=4$, we have no stepwise transitive sets of MOLR, 
and additionally, there are no stepwise transitive $5$-MOLR for $n=7$, $k=5,6$.

Here we note that for each $n\leq7$ all $(n-1)$-tuples of MOLS are stepwise 
transitive, if they exist. Since in each case the maximum set of MOLS for 
these $n$ corresponds to a Galois projective plane (that is, constructed from 
the corresponding Galois field)
this reflects the high degree of symmetry of these planes. For $n=9$ there are 
several projective planes, some of which are not Galois projective planes, and 
we will investigate that case below.

\subsection{Larger orders}
\label{sec_larger}

For $n\geq 8$ we have not generated all $t$-MOLR, even though our programs are 
in principle able to do so. The problem here is that the number of $t$-MOLR becomes 
so large that several peta-byte would be required to store them on disc, and any kind 
of analysis of the whole set would become impractical. Instead, we have focused on 
two interesting subclasses, the stepwise co-isotopic and the stepwise transitive 
$t$-MOLR. These classes are restrictive enough to let us push the generation 
program a few more steps, and we have already seen that they contain a number 
of interesting examples.

In Table~\ref{8x8.} we give the number of stepwise co-isotopic and stepwise
transitive $t$-MOLR for $n=8$ and in Table~\ref{9x9.} we give the number of 
stepwise transitive $t$-MOLR for $n=9$.
For $n=8$ it is clear that for small parameters $k$ and $t$, the stepwise 
co-isotopic $t$-MOLR far outnumber the stepwise transitive ones. We also find  
stepwise co-isotopic $t$-MOLR for all parameters, but not stepwise 
transitive ones. This motivates the following question.

\begin{question}
	For $n\geq7$, is there a stepwise co-isotopic $t$-MOLR for every pair $t,k$ 
	that allows a $t$-MOLR?
\end{question}

\afterpage{
	\clearpage
	\begin{landscape}
		\centering 
		
		\begin{table}[H]
			\begin{center}
				\setlength\tabcolsep{3pt}
				\begin{tabular}{|c||c|c|c|c|c|c|c|}
					\hline 
					&	{2$\times$8}   & {3$\times$8}   & {4$\times$8}   & {5$\times$8}   & {6$\times$8}   & {7$\times$8}   & {8$\times$8}   \\ 
					\hline \hline
					$t=2$	&	 186, 99	&	 \num{446443}, \num{45429}	&	 \num{4432284}, \num{1097655}	&	 \num{3826527}, \num{2569679}	&	 \num{242732}, \num{206612}	&	 484, 305	&	 70, 13	 \\ \hline 
					$t=3$	&	 \num{11565}, 66	&	 \num{9144025}, \num{7627}	&	 \num{178502}, \num{41505}	&	 628, 75	&	 111, 32	&	 10, 6	&	 7, 3	 \\ \hline 
					$t=4$	&	 \num{216950}, 152	&	 \num{1648723}, \num{4284}	&	 \num{3547}, 712	&	 58, 20	&	 4, 0	&	 3, 0	&	 3, 0	 \\ \hline 
					$t=5$	&	 \num{509622}, 19	&	 \num{2652}, 0	&	 267, 0	&	 2, 0	&	 1, 0	&	 1, 0	&	 1, 0	 \\ \hline 
					$t=6$	&	 \num{91013}, 109	&	 975, 908	&	 155, 146	&	 1, 0	&	 1, 0	&	 1, 0	&	 1, 0	 \\ \hline 
					$t=7$	&	 \num{4538}, 5	&	 2, 2	&	 2, 2	&	 1, 1	&	 1, 1	&	 1, 1	&	 1, 1	 \\ \hline 
					
			\end{tabular}\end{center}
			\caption{The number of non-isotopic stepwise co-isotopic and stepwise transitive 
				$t$-MOLR for $n=8$.}
			\label{8x8.}
		\end{table}
		
		\begin{table}[H]
			\begin{center}
				
				\begin{tabular}{|c||c|c|c|c|c|c|c|c|}
					\hline 
					&	{2$\times$9}   & {3$\times$9}   & {4$\times$9}   & {5$\times$9}   & {6$\times$9}   & {7$\times$9}   & {8$\times$9}   & {9$\times$9}   \\ 
					\hline \hline
					$t=2$	&	126	&	\num{1418577}	&	\num{560524587}	&	20019499500	&	67480364637	&	\num{5872237985}	&	\num{14940988}	&	\num{28955}	 \\ \hline 
					$t=3$	&	202	&	\num{72836}	&	\num{1746912}	&	0	&	0	&	0	&	0	&	0	 \\ \hline 
					$t=4$	&	\num{1067}	&	\num{356680}	&	\num{2640163}	&	\num{645453}	&	\num{1816}	&	31	&	7	&	5	 \\ \hline 
					$t=5$	&	17	         &	0	        &	0	                &	0	&	0	&	0	&	0	&	0	 \\ \hline 
					$t=6$	&	543	&	\num{21620}	&	244	       &	33	&	16	&	1	&	1	&	1	 \\ \hline 
					$t=7$	&	39	        &	\num{1532}	&	300	       &	0	&	0	&	0	&	0	&	0	 \\ \hline 
					$t=8$	&	54 	&	48	        &	27	               &	22	&	16	&	9	&	7	&	5	 \\ \hline 
					
			\end{tabular}\end{center}
			\caption{The number of non-isotopic stepwise transitive $t$-MOLR for $n=9$.}
			\label{9x9.}
		\end{table}	
	\end{landscape}
\clearpage
}

For $n=9$ we only have data for the stepwise transitive class, since the number 
of stepwise co-isotopic $t$-MOLR is too large. Here it is clear that the possible 
values of $t$ are quite restricted. We note two interesting facts. 
First, there is a unique stepwise transitive $6$-MOLS, which in turn has unique 
stepwise transitive restrictions to 8 and 7 rows. We present this object in
Figure~\ref{6-9x9}. Second, there are 5 stepwise 
transitive $8$-MOLS, which are presented in Appendix~G.

As mentioned earlier, it is known that there are exactly $4$ 
projective planes of order $9$. The Galois plane corresponds to the $8$-MOLS 
with autotopism group of order \num{10368}, see Table~42 in Appendix~\ref{app_auto9}.   
The other four $8$-MOLS can be divided into two pairs, such that 
both $8$-MOLS in one pair correspond to the Hall plane, and those in the other 
pair correspond to the dual of the Hall plane. In Appendix~\ref{mols9} we 
include data on this pairing. This leaves the 
Hughes plane of order $9$ as the smallest projective plane 
which cannot be defined by a stepwise transitive MOLS. 

With this in mind one may ask about the situation for larger orders as well. 
Wanless~\cite{Wanless_pc} 
has found that $8$ of the $22$ projective planes of order $16$ cannot be 
constructed via a co-isotopic MOLS, and hence they cannot be constructed 
from a stepwise transitive set of MOLS either. In the online catalogue 
provided by Royle (currently available via~\cite{Royle}), these are the 
planes labelled JOHN, BBS4, BBH2 and their duals, BBH1 (which is self-dual), 
and, finally, either MATH or its dual DMATH. In the case of MATH or DMATH, 
the test did not give enough information to discern which one of these two 
planes was constructible in this fashion.

On the other hand, we can prove the following result, where by `the standard 
way', we refer to the construction given in Section~\ref{sec_fingeom}.  Theorem 5.2.5 in~\cite{DK15} shows, 
in our terminology, that the tuples coming from the Galois planes are co-isotopic. However, a closer inspection 
of their proof leads to this stronger statement.
\begin{theorem}\label{thm_galois}
	The set $M$ of $(n-1)$ MOLS $L_1,L_2,\ldots L_{n-1}$ corresponding to a 
	projective plane constructed in the standard way from the finite 
	field $GF(n)$ is stepwise transitive.
\end{theorem}

\begin{proof}	
	Let $a_i = x^{i-1}$ for $i\neq 0$, where $x$ is a generating element 
	of $GF(n)$, and note that column $j$ of $L_k$ has entries	
	\[
		a_0+a_ka_j, a_1+a_ka_j, \ldots, a_{n-1}+a_ka_j
	\]
	in this order, which is identical to column $t$ of $L_1$ for 
	$t=j+k-1$ modulo $n$, so that the columns of $L_k$ are in fact just 
	the columns of $L_1$ permuted cyclically. Rephrasing this, there is an 
	isotopism of $L_1$ that maps it to $L_k$, namely the column permutation $\pi$
	that keeps column $1$ in place, and shifts columns $2, 3, \ldots, n$ of $L_1$ 
	$k-1$ steps forward, modulo $n$. By re-indexing the 
	squares, there is also an isotopism that carries $L_{k_1}$ to $L_{k_2}$ 
	for any $k_1$ and $k_2$. The set $M$ of MOLS is thus co-isotopic.
	
	To see that $M$ is transitive, we need to show that $\pi$ is also 
	an isotopism of the whole set of MOLS $M$. This is clear, though, as 
	$\pi$ leaves column 1 in place, and shifts all the other columns $k-1$ 
	steps forward, modulo $n$. Square $L_{k_1}$ is therefore mapped to square
	$L_{k_1 + k-1}$, where indices are taken modulo $n$.	
	
	When restricting to the first $s$ rows, the argument works in the same 
	way, so the set $M$ is in fact stepwise transitive.
\end{proof}

A full characterization of the family of projective planes which correspond to  
stepwise transitive sets of MOLS would of course be interesting, but even simpler 
questions are left open. 

\begin{question}
	For large $n$, what proportion of the projective planes of order 
	$n$ correspond to co-isotopic or stepwise transitive sets of MOLS? 
\end{question}

Here it seems likely that asymptotically the proportion is 0. It would be of interest 
to use the existing catalogues of finite projective planes, 
the currently most extensive being that of Moorhouse~\cite{Moorhouse} which now 
contains several hundred thousand examples, to check how common these properties are 
among the known non-Galois examples.

\begin{figure}
\begin{tabular}{|ccccccccc|}
\hline 
0 & 1 & 2 & 3 & 4 & 5 & 6 & 7 & 8\\ 
8 & 7 & 6 & 5 & 3 & 2 & 4 & 1 & 0\\ 
7 & 0 & 8 & 4 & 1 & 6 & 5 & 3 & 2\\ 
6 & 2 & 0 & 1 & 8 & 7 & 3 & 5 & 4\\ 
5 & 8 & 4 & 6 & 2 & 3 & 7 & 0 & 1\\ 
4 & 3 & 5 & 7 & 0 & 1 & 8 & 2 & 6\\ 
3 & 4 & 1 & 0 & 7 & 8 & 2 & 6 & 5\\ 
2 & 6 & 7 & 8 & 5 & 0 & 1 & 4 & 3\\ 
1 & 5 & 3 & 2 & 6 & 4 & 0 & 8 & 7\\ 
\hline 
\end{tabular}
\begin{tabular}{|ccccccccc|}
\hline 
0 & 1 & 2 & 3 & 4 & 5 & 6 & 7 & 8\\ 
7 & 8 & 5 & 4 & 2 & 6 & 3 & 0 & 1\\ 
1 & 7 & 6 & 0 & 3 & 8 & 2 & 4 & 5\\ 
2 & 0 & 3 & 6 & 5 & 4 & 1 & 8 & 7\\ 
3 & 4 & 8 & 7 & 1 & 0 & 5 & 6 & 2\\ 
6 & 5 & 1 & 2 & 8 & 7 & 0 & 3 & 4\\ 
8 & 2 & 4 & 5 & 6 & 3 & 7 & 1 & 0\\ 
4 & 3 & 0 & 1 & 7 & 2 & 8 & 5 & 6\\ 
5 & 6 & 7 & 8 & 0 & 1 & 4 & 2 & 3\\ 
\hline 
\end{tabular}
\begin{tabular}{|ccccccccc|}
\hline 
0 & 1 & 2 & 3 & 4 & 5 & 6 & 7 & 8\\ 
6 & 5 & 8 & 0 & 7 & 3 & 1 & 4 & 2\\ 
8 & 4 & 7 & 6 & 2 & 1 & 3 & 5 & 0\\ 
4 & 7 & 5 & 8 & 1 & 2 & 0 & 6 & 3\\ 
7 & 3 & 0 & 2 & 6 & 8 & 4 & 1 & 5\\ 
5 & 6 & 4 & 1 & 3 & 0 & 2 & 8 & 7\\ 
2 & 8 & 6 & 7 & 0 & 4 & 5 & 3 & 1\\ 
1 & 0 & 3 & 5 & 8 & 6 & 7 & 2 & 4\\ 
3 & 2 & 1 & 4 & 5 & 7 & 8 & 0 & 6\\ 
\hline 
\end{tabular}
\begin{tabular}{|ccccccccc|}
\hline 
0 & 1 & 2 & 3 & 4 & 5 & 6 & 7 & 8\\ 
3 & 0 & 7 & 6 & 8 & 1 & 5 & 2 & 4\\ 
6 & 5 & 4 & 8 & 7 & 2 & 0 & 1 & 3\\ 
1 & 6 & 8 & 5 & 0 & 3 & 7 & 4 & 2\\ 
2 & 7 & 1 & 0 & 5 & 4 & 8 & 3 & 6\\ 
7 & 2 & 3 & 4 & 6 & 8 & 1 & 0 & 5\\ 
5 & 3 & 0 & 1 & 2 & 6 & 4 & 8 & 7\\ 
8 & 4 & 5 & 2 & 1 & 7 & 3 & 6 & 0\\ 
4 & 8 & 6 & 7 & 3 & 0 & 2 & 5 & 1\\ 
\hline 
\end{tabular}
\begin{tabular}{|ccccccccc|}
\hline 
0 & 1 & 2 & 3 & 4 & 5 & 6 & 7 & 8\\ 
2 & 6 & 4 & 7 & 1 & 8 & 0 & 3 & 5\\ 
3 & 8 & 0 & 5 & 6 & 7 & 4 & 2 & 1\\ 
5 & 4 & 7 & 2 & 3 & 0 & 8 & 1 & 6\\ 
8 & 5 & 3 & 1 & 7 & 6 & 2 & 4 & 0\\ 
1 & 0 & 6 & 8 & 2 & 4 & 7 & 5 & 3\\ 
4 & 7 & 8 & 6 & 5 & 1 & 3 & 0 & 2\\ 
6 & 2 & 1 & 4 & 0 & 3 & 5 & 8 & 7\\ 
7 & 3 & 5 & 0 & 8 & 2 & 1 & 6 & 4\\ 
\hline 
\end{tabular}\hfill
\begin{tabular}{|ccccccccc|}
\hline 
0 & 1 & 2 & 3 & 4 & 5 & 6 & 7 & 8\\ 
1 & 4 & 0 & 8 & 6 & 7 & 2 & 5 & 3\\ 
2 & 6 & 5 & 7 & 8 & 3 & 1 & 0 & 4\\ 
8 & 3 & 1 & 4 & 7 & 6 & 5 & 2 & 0\\ 
4 & 2 & 6 & 5 & 0 & 1 & 3 & 8 & 7\\ 
3 & 8 & 7 & 0 & 5 & 2 & 4 & 6 & 1\\ 
7 & 5 & 3 & 2 & 1 & 0 & 8 & 4 & 6\\ 
5 & 7 & 8 & 6 & 3 & 4 & 0 & 1 & 2\\ 
6 & 0 & 4 & 1 & 2 & 8 & 7 & 3 & 5\\ 	
\hline 
\end{tabular} 
\caption{The unique stepwise transitive $6$-MOLS of order 9.}
\label{6-9x9}
\end{figure}

\subsection{Autotopism group sizes of $\mathbf{t}$-MOLR}
\label{sec_auto}

We have computed the order of the autotopism group for all sets of MOLR 
discussed so far in the paper. Detailed statistics of these orders are 
given in Appendices A to F. We will here discuss some of the symmetry 
properties of sets of MOLR in general, and some additional observations 
based on our data.

First let us note that the case of $2\times n$ sets of MOLR is somewhat special. 
If we follow the construction for a partial net using a $2\times n$ $t$-MOLR 
$A$, each of the lines which do not correspond to a row has two points, and no 
such line connects two vertices in the same row. This means that if we delete 
the two lines with $n$ points, we have a bipartite graph $g(A)$, where the two rows 
give us the bipartition. Additionally, the edges coming from each rectangle add 
a perfect matching, as do the edges coming from the columns, so we have a 
$(t+1)$-regular bipartite graph with a natural edge colouring given by these 
matchings. The autotopism group of $A$ now corresponds to automorphisms of this 
edge-coloured graph which map the matching given by the column-lines to itself.
If we assume that $A$ is normalized, we can also invert this construction and 
reconstruct the $t$-MOLR $A$. Now, for $t=n-1$ this defines a proper edge 
colouring of $K_{n,n}$, i.e., a Latin square, so here we obtain a mapping from 
$2\times n$ $(n-1)$-MOLR to a Latin square $L(A)$, and an autotopism of $A$ 
defines an autotopism of $L(A)$ which fixes one symbol, again corresponding 
to the column lines in the partial net.

Given a regularity property we may also look at how it interplays with 
restrictions of a set of MOLR. First let us note that any subset of rectangles 
from a co-isotopic $t$-MOLR is co-isotopic, so this case is trivial, and the 
same is true for a stepwise co-isotopic $t$-MOLR.  For transitivity, the group 
structure comes into play. Given a transitive $t$-MOLR $A$ with autotopism group 
$G$ we get a subgroup $G'$ which describes the action of $G$ on the set of rectangles. 
If $G'$ has an element $g$ of order $r$ we will obtain a transitive $r$-MOLR from $A$ by 
taking the orbit of a single rectangle from $A$ under $g$. Whenever we have 
a transitive $t$-MOLR with autotopism group of order $t$ this implies the 
existence of $p$-MOLR for the same size $k\times n$ for every prime factor $p$ 
of $t$. For general parameters we observe that when $t$ or $k$ is small, transitive 
$t$-MOLR with autotopism group of order exactly $t$ are common.

Following this, we say that a transitive $t$-MOLR $A$ is $G$-complete if there does 
not exist a $t'>t$ and a $t'$-MOLR $B$, with $A\subsetneq B$, such that $A$ is the orbit of a 
rectangle in $B$ under an element $g\in \mathrm{Aut}(B)$, and otherwise we say 
that $A$ is $G$-incomplete. As noted, there are many examples of $t$-MOLR with 
autotopism groups of size $t$, and hence we will also have many incomplete 
$r$-MOLR of the same size and $r$ a divisor of $t$. 

\begin{observation} In our data we found the following:
	\begin{enumerate}
		\item Among the stepwise transitive sets of MOLR for $n=8$ there is one 
			$3$-MOLS with autotopism group of order 48, which is $G$-complete. We 
			display that example in Figure~\ref{Gcomp83}.  
		
		\item For $n=9$ none of the stepwise transitive $4$-MOLR are $G$-complete.   
			The $4$-MOLR with autotopism groups of orders \num{5184} and \num{2592} 
			both correspond to the $8$-MOLS with autotopism group of order 
			\num{10368}. 
			The two $4$-MOLR with autotopism group of order 64 correspond to 
			the two $8$-MOLR with autotopism group of order 384.
	
		\item  For $n=9$ the $8$-MOLS with autotopism group of order 
			\num{31104} does not correspond to any $G$-incomplete stepwise 
			transitive $4$-MOLS.
		
		\item For $n=9$ the stepwise transitive set of 6 MOLS is $G$-complete.
	\end{enumerate}
\end{observation}

\begin{figure}
\begin{tabular}{|c|c|c|c|c|c|c|c|}
	\hline 
	0 & 1 & 2 & 3 & 4 & 5 & 6 & 7  \\ \hline 
	7 & 6 & 5 & 4 & 3 & 2 & 1 & 0  \\ \hline 
	6 & 7 & 4 & 5 & 2 & 3 & 0 & 1  \\ \hline 
	5 & 4 & 7 & 6 & 1 & 0 & 3 & 2  \\ \hline 
	4 & 5 & 6 & 7 & 0 & 1 & 2 & 3  \\ \hline 
	3 & 2 & 1 & 0 & 7 & 6 & 5 & 4  \\ \hline 
	2 & 0 & 3 & 1 & 6 & 4 & 7 & 5  \\ \hline 
	1 & 3 & 0 & 2 & 5 & 7 & 4 & 6  \\ \hline 
\end{tabular}
\begin{tabular}{|c|c|c|c|c|c|c|c|}
\hline 
 0 & 1 & 2 & 3 & 4 & 5 & 6 & 7   \\ \hline 
 6 & 7 & 4 & 5 & 2 & 3 & 0 & 1   \\ \hline 
 5 & 4 & 7 & 6 & 1 & 0 & 3 & 2   \\ \hline 
 2 & 3 & 0 & 1 & 6 & 7 & 4 & 5   \\ \hline 
 1 & 0 & 3 & 2 & 5 & 4 & 7 & 6   \\ \hline 
 7 & 6 & 5 & 4 & 3 & 2 & 1 & 0   \\ \hline 
 4 & 2 & 1 & 7 & 0 & 6 & 5 & 3   \\ \hline 
 3 & 5 & 6 & 0 & 7 & 1 & 2 & 4   \\ \hline 
\end{tabular}
\begin{tabular}{|c|c|c|c|c|c|c|c|}
	\hline 
 0 & 1 & 2 & 3 & 4 & 5 & 6 & 7  \\ \hline 
 1 & 0 & 3 & 2 & 5 & 4 & 7 & 6  \\ \hline 
 3 & 2 & 1 & 0 & 7 & 6 & 5 & 4  \\ \hline 
 7 & 6 & 5 & 4 & 3 & 2 & 1 & 0  \\ \hline 
 5 & 4 & 7 & 6 & 1 & 0 & 3 & 2  \\ \hline 
 4 & 5 & 6 & 7 & 0 & 1 & 2 & 3  \\ \hline 
 6 & 3 & 0 & 5 & 2 & 7 & 4 & 1  \\ \hline 
 2 & 7 & 4 & 1 & 6 & 3 & 0 & 5  \\ \hline 
\end{tabular} 
\caption{The $G$-complete $3$-MOLS with $n=8$}
\label{Gcomp83}
\end{figure}

For stepwise transitive sets of MOLR, restrictions become far less well-behaved.
Given a transitive $t$-MOLR $A$ and an autotopism $g$ which has order 
$r$ on the set of rectangles, we know that we will obtain a transitive 
$r$-MOLR $A'$ from $A$. However, assuming that $A$ is stepwise transitive 
does not necessarily lead to stepwise transitivity for $A'$. In order 
for this to happen it must also be the case that each of the stepwise 
transitive sets of MOLR which are used to construct $A$ have autotopisms with 
the same orbit as $g$, and this is not always the case. We see one such 
example at $n=9$, where a stepwise transitive $6$-MOLS exists, 
but no stepwise transitive triple.

\section{Concluding remarks}
\label{sec_concl}

In this paper we have focused on enumeration of $t$-MOLR up to $n=9$, 
coming tantalisingly close to the, in this setting, special value $n=10$. 
As we have mentioned, a significant theoretical and computational effort 
led to the result that there is no set of 9 mutually orthogonal Latin 
squares of order 10, and hence no projective plane of order 10.   
However, an even more basic question remains:  

\begin{question}\label{Q_triples}
	Is there a triple of mutually orthogonal Latin squares of order 10?
\end{question}

A large number of pairwise orthogonal Latin squares of order 10 are known, 
and as mentioned above, the self-orthogonal Latin squares of order 10 have 
been completely enumerated. An exhaustive search by McKay, Meynert, and 
Myrvold~\cite{MMM} proved that no square of order 10 with non-trivial 
autoparatopism group is part of an orthogonal triple. The total number of 
Latin squares of order 10 with trivial autoparatopism group was however too 
large for a complete search for orthogonal triples.

Note that these results do not immediately exclude the 
existence of transitive, or even stepwise transitive, triples of MOLS or 
MOLR of order 10, since the autotopism group of a single square or rectangle 
in such a triple can be trivial. Our data provides several such examples for 
rectangles.  However, as an anonymous reviewer pointed out, from a transitive 
triple of MOLS one can, by translation to an orthogonal array, construct a 
triple of MOLS where at least one square has non-trivial autoparatopisms.  
Thus the result of McKay, Meynert, and Myrvold~\cite{MMM} does in fact 
rule out the existence of transitive triples of MOLS of order 10 as well.

There are other restricted versions of Question~\ref{Q_triples} which remain.

\begin{question}\label{q_co-isotopic}
	Is there a co-isotopic $3$-MOLS of order $10$?  
\end{question} 

\begin{question}\label{q_stepwise}
	Is there a stepwise co-isotopic $3$-MOLS of order $10$? 
\end{question} 

A negative answer to Question~\ref{q_co-isotopic} would lead to another 
extension of McKay, Meynert, and Myrvold's result from~\cite{MMM}. For 
Question~\ref{q_stepwise}, a 
more specialised version of the type of search we have performed might be 
able to handle the case $t=3$ for $n=10$ as well.

In~\cite{WanlessEgan} Egan and Wanless tried to find an example of $3$ Latin 
squares of order $10$ that come as close to being mutually orthogonal as possible. 
They presented an example of $3$ Latin squares such that the first is orthogonal 
to the other two, and the final two produce $91$ different symbol pairs when
superimposed. There have also been earlier examples 
of MOLR and pairwise almost orthogonal Latin squares for $n=10$. In~\cite{franklin84}, Franklin constructed
examples of triples of pairwise orthogonal $9\times 10$ rectangles, 
in order to be used in the construction of designs, and in~\cite{WANLESS2001} 
Wanless constructed a set of $4$ such rectangles.

Using our program we performed a partial search for stepwise transitive MOLR 
with $n=10$. We 
found several examples of stepwise transitive triples of $8\times 10$ MOLR, 
and some of these could be extended to triples of $9\times 10$ MOLR, but not 
while preserving transitivity. In Figure~\ref{molr10} we give one such 
example. This example can be uniquely extended to $3$ Latin squares, such 
that all positions which break orthogonality lie in the last row. 
Unfortunately, none of the examples we found could be extended to a triple 
of MOLS.    

\begin{figure}
	\begin{tabular}{|c|c|c|c|c|c|c|c|c|c|}
		\hline 
		0 & 1 & 2 & 3 & 4 & 5 & 6 & 7 & 8 & 9 \\ \hline 
		9 & 8 & 7 & 6 & 5 & 4 & 3 & 1 & 0 & 2 \\ \hline 
		8 & 9 & 3 & 5 & 7 & 1 & 0 & 2 & 6 & 4 \\ \hline 
		7 & 6 & 4 & 0 & 3 & 2 & 9 & 8 & 5 & 1 \\ \hline 
		6 & 7 & 5 & 4 & 1 & 8 & 2 & 0 & 9 & 3 \\ \hline 
		5 & 4 & 6 & 1 & 9 & 0 & 8 & 3 & 2 & 7 \\ \hline 
		4 & 5 & 0 & 7 & 2 & 6 & 1 & 9 & 3 & 8 \\ \hline 
		3 & 0 & 1 & 2 & 8 & 9 & 4 & 6 & 7 & 5 \\ \hline 
		2 & 3 & 8 & 9 & 6 & 7 & 5 & 4 & 1 & 0 \\ \hline 
	\end{tabular}	
	\begin{tabular}{|c|c|c|c|c|c|c|c|c|c|}
		\hline 
		0 & 1 & 2 & 3 & 4 & 5 & 6 & 7 & 8 & 9 \\ \hline 
		8 & 9 & 6 & 5 & 0 & 7 & 2 & 4 & 3 & 1 \\ \hline 
		4 & 7 & 0 & 6 & 2 & 3 & 1 & 8 & 9 & 5 \\ \hline 
		5 & 3 & 8 & 9 & 1 & 4 & 0 & 2 & 7 & 6 \\ \hline 
		2 & 4 & 3 & 1 & 8 & 0 & 9 & 5 & 6 & 7 \\ \hline 
		1 & 0 & 4 & 2 & 3 & 6 & 7 & 9 & 5 & 8 \\ \hline 
		9 & 2 & 7 & 0 & 6 & 8 & 5 & 1 & 4 & 3 \\ \hline 
		6 & 8 & 9 & 7 & 5 & 2 & 3 & 0 & 1 & 4 \\ \hline 
		3 & 5 & 1 & 4 & 7 & 9 & 8 & 6 & 0 & 2 \\ \hline 
	\end{tabular}	
	\begin{tabular}{|c|c|c|c|c|c|c|c|c|c|}
		\hline 
		0 & 1 & 2 & 3 & 4 & 5 & 6 & 7 & 8 & 9\\ \hline 
		4 & 5 & 3 & 2 & 6 & 9 & 1 & 8 & 7 & 0\\ \hline 
		9 & 6 & 7 & 0 & 8 & 4 & 2 & 5 & 1 & 3\\ \hline 
		1 & 9 & 0 & 4 & 5 & 3 & 8 & 6 & 2 & 7\\ \hline 
		3 & 0 & 8 & 6 & 2 & 1 & 7 & 9 & 5 & 4\\ \hline 
		7 & 2 & 5 & 9 & 1 & 8 & 3 & 0 & 4 & 6\\ \hline 
		8 & 4 & 1 & 5 & 9 & 7 & 0 & 3 & 6 & 2\\ \hline 
		2 & 3 & 6 & 8 & 7 & 0 & 5 & 4 & 9 & 1\\ \hline 
		6 & 8 & 4 & 7 & 0 & 2 & 9 & 1 & 3 & 5\\ \hline 
	\end{tabular}
\caption{A $3$-MOLR of size $9\times 10$, whose restriction to the first
	8 rows is stepwise transitive.}\label{molr10}
\end{figure}

\section*{Acknowledgments}

The authors would like to thank Ian Wanless for providing us with useful 
references and data on projective planes. We would also like to thank one 
anonymous reviewer for pointing out the connection to orthogonal arrays 
mentioned in the last section, and the other diligent reviewers for many 
helpful comments.

The computational work was performed on resources provided by the Swedish 
National Infrastructure for Computing (SNIC) at 
High Performance Computing Center North (HPC2N). 
This work was supported by the Swedish strategic research programme eSSENCE.    
This work was supported by The Swedish Research Council grant \num{2014}-\num{4897}. 

\bibliographystyle{plain}

\newpage


\appendix

\section{Sizes of autotopism groups of MOLR for $\mathbf{n=4}$}
\label{app_auto4}

\begin{minipage}{.5\linewidth}
\begin{table}[H]

\begin{tabular}{|c|c||c|c|}
			\hline 
			2$\times$4 & $\Sigma$ & 8 & 16 \\ \hline 
			non-ISO      & 3 & 1 & 2 \\ \hline 
			co-iso.      & 2 & 0 & 2 \\ \hline 
			transitive      & 2 & 0 & 2 \\ \hline 
			st. co-iso.      & 2 & 0 & 2 \\ \hline 
			st. trans.      & 2 & 0 & 2 \\ \hline 
		\end{tabular}
		
		\begin{tabular}{|c|c||c|c|}
			\hline 
			3$\times$4 & $\Sigma$ & 8 & 24 \\ \hline 
			non-ISO      & 2 & 1 & 1 \\ \hline 
			co-iso.      & 2 & 1 & 1 \\ \hline 
			transitive      & 2 & 1 & 1 \\ \hline 
			st. co-iso.      & 1 & 0 & 1 \\ \hline 
			st. trans.      & 1 & 0 & 1 \\ \hline 
		\end{tabular}
		
		\begin{tabular}{|c|c||c|}
			\hline 
			4$\times$4 & $\Sigma$ & 96 \\ \hline 
			non-ISO      & 1 & 1 \\ \hline 
			co-iso.      & 1 & 1 \\ \hline 
			transitive      & 1 & 1 \\ \hline 
			st. co-iso.      & 1 & 1 \\ \hline 
			st. trans.      & 1 & 1 \\ \hline 
		\end{tabular}
		
	\caption{$2$-MOLR for $n=4$.}
	\label{4x2.}
\end{table}
\end{minipage}
\begin{minipage}{.5\linewidth}
\begin{table}[H]
	\begin{center}
		
		\begin{tabular}{|c|c||c|c|}
			\hline 
			2$\times$4 & $\Sigma$ & 16 & 48 \\ \hline 
			non-ISO      & 2 & 1 & 1 \\ \hline 
			co-iso.      & 1 & 0 & 1 \\ \hline 
			transitive      & 1 & 0 & 1 \\ \hline 
			st. co-iso.      & 1 & 0 & 1 \\ \hline 
			st. trans.      & 1 & 0 & 1 \\ \hline 
		\end{tabular}
		
		\begin{tabular}{|c|c||c|}
			\hline 
			3$\times$4 & $\Sigma$ & 72 \\ \hline 
			non-ISO      & 1 & 1 \\ \hline 
			co-iso.      & 1 & 1 \\ \hline 
			transitive      & 1 & 1 \\ \hline 
			st. co-iso.      & 1 & 1 \\ \hline 
			st. trans.      & 1 & 1 \\ \hline 
		\end{tabular}
		
		\begin{tabular}{|c|c||c|}
			\hline 
			4$\times$4 & $\Sigma$ & 288 \\ \hline 
			non-ISO      & 1 & 1 \\ \hline 
			co-iso.      & 1 & 1 \\ \hline 
			transitive      & 1 & 1 \\ \hline 
			st. co-iso.      & 1 & 1 \\ \hline 
			st. trans.      & 1 & 1 \\ \hline 
		\end{tabular}
		
	\end{center}
	\caption{$3$-MOLR for $n=4$.}
	\label{4x3.}
\end{table}
\end{minipage}

\section{Sizes of autotopism groups of MOLR for $\mathbf{n=5}$}
\label{app_auto5}
	
\begin{table}[H]
		
		\begin{tabular}{|c|c||c|c|c|c|}
			\hline 
			2$\times$5 & $\Sigma$ & 2 & 4 & 10 & 20 \\ \hline 
			non-ISO      & 5 & 1 & 2 & 1 & 1 \\ \hline 
			co-iso.      & 4 & 0 & 2 & 1 & 1 \\ \hline 
			transitive      & 3 & 0 & 2 & 0 & 1 \\ \hline 
			st. co-iso.      & 4 & 0 & 2 & 1 & 1 \\ \hline 
			st. trans.      & 3 & 0 & 2 & 0 & 1 \\ \hline 
		\end{tabular}
		\vfill
		\begin{tabular}{|c|c||c|c|c|c|c|}
			\hline 
			3$\times$5 & $\Sigma$ & 1 & 2 & 4 & 10 & 20 \\ \hline 
			non-ISO      & 14 & 1 & 7 & 3 & 2 & 1 \\ \hline 
			co-iso.      & 11 & 1 & 4 & 3 & 2 & 1 \\ \hline 
			transitive      & 9 & 0 & 4 & 3 & 1 & 1 \\ \hline 
			st. co-iso.      & 7 & 0 & 1 & 3 & 2 & 1 \\ \hline 
			st. trans.      & 6 & 0 & 1 & 3 & 1 & 1 \\ \hline 
		\end{tabular}
		
		\begin{tabular}{|c|c||c|c|}
			\hline 
			4$\times$5 & $\Sigma$ & 20 & 40 \\ \hline 
			non-ISO      & 2 & 1 & 1 \\ \hline 
			co-iso.      & 2 & 1 & 1 \\ \hline 
			transitive      & 1 & 0 & 1 \\ \hline 
			st. co-iso.      & 2 & 1 & 1 \\ \hline 
			st. trans.      & 1 & 0 & 1 \\ \hline 
		\end{tabular}
		
		\begin{tabular}{|c|c||c|c|}
			\hline 
			5$\times$5 & $\Sigma$ & 100 & 200 \\ \hline 
			non-ISO      & 2 & 1 & 1 \\ \hline 
			co-iso.      & 2 & 1 & 1 \\ \hline 
			transitive      & 1 & 0 & 1 \\ \hline 
			st. co-iso.      & 2 & 1 & 1 \\ \hline 
			st. trans.      & 1 & 0 & 1 \\ \hline 
		\end{tabular}
		
	\caption{$2$-MOLR for $n=5$.}
	\label{5x2.}
\end{table}
\vfill
\begin{table}[H]
		
		\begin{tabular}{|c|c||c|c|c|}
			\hline 
			2$\times$5 & $\Sigma$ & 2 & 6 & 10 \\ \hline 
			non-ISO      & 4 & 1 & 2 & 1 \\ \hline 
			co-iso.      & 3 & 0 & 2 & 1 \\ \hline 
			transitive      & 2 & 0 & 2 & 0 \\ \hline 
			st. co-iso.      & 3 & 0 & 2 & 1 \\ \hline 
			st. trans.      & 2 & 0 & 2 & 0 \\ \hline 
		\end{tabular}
		
		\begin{tabular}{|c|c||c|}
			\hline 
			3$\times$5 & $\Sigma$ & 10 \\ \hline 
			non-ISO      & 1 & 1 \\ \hline 
			co-iso.      & 1 & 1 \\ \hline 
			st. co-iso.      & 1 & 1 \\ \hline 
		\end{tabular}
		
		\begin{tabular}{|c|c||c|}
			\hline 
			4$\times$5 & $\Sigma$ & 20 \\ \hline 
			non-ISO      & 1 & 1 \\ \hline 
			co-iso.      & 1 & 1 \\ \hline 
			st. co-iso.      & 1 & 1 \\ \hline 
		\end{tabular}
		
		\begin{tabular}{|c|c||c|}
			\hline 
			5$\times$5 & $\Sigma$ & 100 \\ \hline 
			non-ISO      & 1 & 1 \\ \hline 
			co-iso.      & 1 & 1 \\ \hline 
		\end{tabular}
		
	\caption{$3$-MOLR for $n=5$.}
	\label{5x3.}
\end{table}

\begin{table}[H]
		
		\begin{tabular}{|c|c||c|c|c|}
			\hline 
			2$\times$5 & $\Sigma$ & 6 & 24 & 40 \\ \hline 
			non-ISO      & 3 & 1 & 1 & 1 \\ \hline 
			co-iso.      & 2 & 0 & 1 & 1 \\ \hline 
			transitive      & 2 & 0 & 1 & 1 \\ \hline 
			st. co-iso.      & 2 & 0 & 1 & 1 \\ \hline 
			st. trans.      & 2 & 0 & 1 & 1 \\ \hline 
		\end{tabular}
		
		\begin{tabular}{|c|c||c|}
			\hline 
			3$\times$5 & $\Sigma$ & 40 \\ \hline 
			non-ISO      & 1 & 1 \\ \hline 
			co-iso.      & 1 & 1 \\ \hline 
			transitive      & 1 & 1 \\ \hline 
			st. co-iso.      & 1 & 1 \\ \hline 
			st. trans.      & 1 & 1 \\ \hline 
		\end{tabular}
		
		\begin{tabular}{|c|c||c|}
			\hline 
			4$\times$5 & $\Sigma$ & 80 \\ \hline 
			non-ISO      & 1 & 1 \\ \hline 
			co-iso.      & 1 & 1 \\ \hline 
			transitive      & 1 & 1 \\ \hline 
			st. co-iso.      & 1 & 1 \\ \hline 
			st. trans.      & 1 & 1 \\ \hline 
		\end{tabular}
		
		\begin{tabular}{|c|c||c|}
			\hline 
			5$\times$5 & $\Sigma$ & 400 \\ \hline 
			non-ISO      & 1 & 1 \\ \hline 
			co-iso.      & 1 & 1 \\ \hline 
			transitive      & 1 & 1 \\ \hline 
			st. co-iso.      & 1 & 1 \\ \hline 
			st. trans.      & 1 & 1 \\ \hline 
		\end{tabular}
		
	\caption{$4$-MOLR for $n=5$.}
	\label{5x4.}
\end{table}




\section{Sizes of autotopism groups of MOLR for $\mathbf{n=6}$}	
\label{app_auto6}

\begin{table}[H]
		\setlength\tabcolsep{3pt}
		\begin{tabular}{|c|c||c|c|c|c|c|c|c|c|}
			\hline 
			2$\times$6 & $\Sigma$ & 1 & 2 & 4 & 6 & 8 & 12 & 24 & 72 \\ \hline 
			non-ISO      & 28 & 1 & 7 & 7 & 1 & 3 & 6 & 2 & 1 \\ \hline 
			co-iso.      & 12 & 0 & 3 & 2 & 0 & 3 & 1 & 2 & 1 \\ \hline 
			transitive      & 11 & 0 & 2 & 2 & 0 & 3 & 1 & 2 & 1 \\ \hline 
			st. co-iso.      & 12 & 0 & 3 & 2 & 0 & 3 & 1 & 2 & 1 \\ \hline 
			st. trans.      & 11 & 0 & 2 & 2 & 0 & 3 & 1 & 2 & 1 \\ \hline 
		\end{tabular}
		
		\begin{tabular}{|c|c||c|c|c|c|c|c|c|c|c|}
			\hline 
			3$\times$6 & $\Sigma$ & 1 & 2 & 3 & 4 & 6 & 8 & 12 & 24 & 216 \\ \hline 
			non-ISO      & \num{1526} & \num{1155} & 252 & 18 & 59 & 19 & 8 & 11 & 3 & 1 \\ \hline 
			co-iso.      & 280 & 89 & 117 & 1 & 40 & 13 & 8 & 8 & 3 & 1 \\ \hline 
			transitive      & 170 & 0 & 100 & 0 & 40 & 10 & 8 & 8 & 3 & 1 \\ \hline 
			st. co-iso.      & 158 & 43 & 63 & 1 & 26 & 9 & 8 & 4 & 3 & 1 \\ \hline 
			st. trans.      & 103 & 0 & 52 & 0 & 26 & 9 & 8 & 4 & 3 & 1 \\ \hline 
		\end{tabular}
		
		\begin{tabular}{|c|c||c|c|c|c|c|c|c|c|c|c|}
			\hline 
			4$\times$6 & $\Sigma$ & 1 & 2 & 3 & 4 & 6 & 8 & 12 & 16 & 24 & 48 \\ \hline 
			non-ISO      & \num{2036} & \num{1425} & 425 & 30 & 78 & 35 & 16 & 21 & 1 & 3 & 2 \\ \hline 
			co-iso.      & 229 & 36 & 112 & 5 & 31 & 11 & 16 & 12 & 1 & 3 & 2 \\ \hline 
			transitive      & 160 & 0 & 92 & 0 & 29 & 9 & 15 & 9 & 1 & 3 & 2 \\ \hline 
			st. co-iso.      & 66 & 7 & 27 & 1 & 11 & 3 & 8 & 6 & 0 & 2 & 1 \\ \hline 
			st. trans.      & 50 & 0 & 22 & 0 & 10 & 3 & 6 & 6 & 0 & 2 & 1 \\ \hline 
		\end{tabular}
		
		\begin{tabular}{|c|c||c|c|c|c|c|c|c|c|c|c|}
			\hline 
			5$\times$6 & $\Sigma$ & 1 & 2 & 3 & 4 & 6 & 8 & 12 & 16 & 24 & 48 \\ \hline 
			non-ISO      & 85 & 5 & 25 & 2 & 26 & 4 & 11 & 7 & 2 & 1 & 2 \\ \hline 
			co-iso.      & 43 & 0 & 10 & 0 & 13 & 3 & 6 & 7 & 2 & 0 & 2 \\ \hline 
			transitive      & 36 & 0 & 6 & 0 & 11 & 2 & 6 & 7 & 2 & 0 & 2 \\ \hline 
			st. co-iso.      & 13 & 0 & 2 & 0 & 0 & 1 & 4 & 2 & 2 & 0 & 2 \\ \hline 
			st. trans.      & 12 & 0 & 1 & 0 & 0 & 1 & 4 & 2 & 2 & 0 & 2 \\ \hline 
		\end{tabular}
		
	\caption{$2$-MOLR for $n=6$.}
	\label{6x2.}
\end{table}

\begin{table}[H]
		\setlength\tabcolsep{3pt}
		\begin{tabular}{|c|c||c|c|c|c|c|c|c|c|c|}
			\hline 
			2$\times$6 & $\Sigma$ & 1 & 2 & 4 & 6 & 8 & 12 & 24 & 36 & 72 \\ \hline 
			non-ISO      & 103 & 24 & 25 & 26 & 2 & 7 & 13 & 4 & 1 & 1 \\ \hline 
			co-iso.      & 16 & 3 & 4 & 3 & 0 & 0 & 3 & 1 & 1 & 1 \\ \hline 
			transitive      & 6 & 0 & 0 & 0 & 0 & 0 & 3 & 1 & 1 & 1 \\ \hline 
			st. co-iso.      & 16 & 3 & 4 & 3 & 0 & 0 & 3 & 1 & 1 & 1 \\ \hline 
			st. trans.      & 6 & 0 & 0 & 0 & 0 & 0 & 3 & 1 & 1 & 1 \\ \hline 
		\end{tabular}
		
		\begin{tabular}{|c|c||c|c|c|c|c|c|c|c|}
			\hline 
			3$\times$6 & $\Sigma$ & 1 & 2 & 3 & 4 & 6 & 12 & 18 & 36 \\ \hline 
			non-ISO      & \num{2572} & \num{1980} & 442 & 54 & 27 & 55 & 6 & 4 & 4 \\ \hline 
			co-iso.      & 115 & 41 & 32 & 11 & 2 & 18 & 3 & 4 & 4 \\ \hline 
			transitive      & 29 & 0 & 0 & 6 & 0 & 13 & 2 & 4 & 4 \\ \hline 
			st. co-iso.      & 32 & 11 & 11 & 2 & 2 & 2 & 1 & 1 & 2 \\ \hline 
			st. trans.      & 4 & 0 & 0 & 0 & 0 & 0 & 1 & 1 & 2 \\ \hline 
		\end{tabular}
		
		\begin{tabular}{|c|c||c|c|c|c|c|c|c|c|c|c|c|}
			\hline 
			4$\times$6 & $\Sigma$ & 1 & 2 & 3 & 4 & 6 & 8 & 9 & 12 & 18 & 24 & 36 \\ \hline 
			non-ISO      & 513 & 93 & 194 & 96 & 37 & 64 & 3 & 2 & 11 & 9 & 1 & 3 \\ \hline 
			co-iso.      & 62 & 1 & 8 & 11 & 1 & 23 & 0 & 2 & 3 & 9 & 1 & 3 \\ \hline 
			transitive      & 39 & 0 & 0 & 3 & 0 & 18 & 0 & 2 & 3 & 9 & 1 & 3 \\ \hline 
			st. co-iso.      & 4 & 0 & 0 & 0 & 0 & 3 & 0 & 0 & 0 & 0 & 0 & 1 \\ \hline 
			st. trans.      & 1 & 0 & 0 & 0 & 0 & 0 & 0 & 0 & 0 & 0 & 0 & 1 \\ \hline 
		\end{tabular}
		
		\begin{tabular}{|c|c||c|c|c|c|}
			\hline 
			5$\times$6 & $\Sigma$ & 3 & 6 & 9 & 18 \\ \hline 
			non-ISO      & 7 & 2 & 2 & 1 & 2 \\ \hline 
			co-iso.      & 4 & 0 & 1 & 1 & 2 \\ \hline 
			transitive      & 3 & 0 & 0 & 1 & 2 \\ \hline 
		\end{tabular}
		
	\caption{$3$-MOLR for $n=6$.}
	\label{6x3.}
\end{table}

\begin{table}[H]
		\setlength\tabcolsep{3pt}
		\begin{tabular}{|c|c||c|c|c|c|c|c|c|c|c|c|c|}
			\hline 
			2$\times$6 & $\Sigma$ & 1 & 2 & 3 & 4 & 8 & 12 & 16 & 18 & 24 & 36 & 48 \\ \hline 
			non-ISO      & 92 & 14 & 18 & 1 & 28 & 8 & 10 & 2 & 1 & 4 & 4 & 2 \\ \hline 
			co-iso.      & 9 & 0 & 1 & 0 & 2 & 2 & 0 & 2 & 0 & 0 & 0 & 2 \\ \hline 
			transitive      & 8 & 0 & 0 & 0 & 2 & 2 & 0 & 2 & 0 & 0 & 0 & 2 \\ \hline 
			st. co-iso.      & 9 & 0 & 1 & 0 & 2 & 2 & 0 & 2 & 0 & 0 & 0 & 2 \\ \hline 
			st. trans.      & 8 & 0 & 0 & 0 & 2 & 2 & 0 & 2 & 0 & 0 & 0 & 2 \\ \hline 
		\end{tabular}
		
		\begin{tabular}{|c|c||c|c|c|c|c|c|c|c|c|c|c|c|c|}
			\hline 
			3$\times$6 & $\Sigma$ & 1 & 2 & 3 & 4 & 6 & 8 & 12 & 16 & 18 & 24 & 36 & 48 & 144 \\ \hline 
			non-ISO      & 118 & 16 & 38 & 5 & 22 & 10 & 6 & 9 & 2 & 4 & 1 & 3 & 1 & 1 \\ \hline 
			co-iso.      & 19 & 0 & 0 & 0 & 5 & 1 & 6 & 2 & 2 & 0 & 1 & 0 & 1 & 1 \\ \hline 
			transitive      & 17 & 0 & 0 & 0 & 5 & 0 & 6 & 1 & 2 & 0 & 1 & 0 & 1 & 1 \\ \hline 
			st. co-iso.      & 15 & 0 & 0 & 0 & 3 & 0 & 6 & 1 & 2 & 0 & 1 & 0 & 1 & 1 \\ \hline 
			st. trans.      & 15 & 0 & 0 & 0 & 3 & 0 & 6 & 1 & 2 & 0 & 1 & 0 & 1 & 1 \\ \hline 
		\end{tabular}
		
		\begin{tabular}{|c|c||c|c|c|c|c|c|c|}
			\hline 
			4$\times$6 & $\Sigma$ & 3 & 6 & 9 & 12 & 18 & 24 & 72 \\ \hline 
			non-ISO      & 12 & 2 & 2 & 2 & 1 & 3 & 1 & 1 \\ \hline 
			co-iso.      & 4 & 0 & 0 & 1 & 1 & 0 & 1 & 1 \\ \hline 
			transitive      & 3 & 0 & 0 & 0 & 1 & 0 & 1 & 1 \\ \hline 
		\end{tabular}
		
		\begin{tabular}{|c|c||c|c|c|c|c|c|}
			\hline 
			5$\times$6 & $\Sigma$ & 6 & 9 & 18 & 24 & 36 & 72 \\ \hline 
			non-ISO      & 8 & 1 & 2 & 1 & 1 & 1 & 2 \\ \hline 
			co-iso.      & 4 & 0 & 0 & 0 & 1 & 1 & 2 \\ \hline 
			transitive      & 4 & 0 & 0 & 0 & 1 & 1 & 2 \\ \hline 
		\end{tabular}
		
	\caption{$4$-MOLR for $n=6$.}
	\label{6x4.}
\end{table}

\begin{table}[H]
		\setlength\tabcolsep{3pt}
		\begin{tabular}{|c|c||c|c|c|c|c|c|c|c|c|c|c|c|c|}
			\hline 
			2$\times$6 & $\Sigma$ & 1 & 2 & 4 & 6 & 8 & 12 & 16 & 20 & 24 & 36 & 48 & 72 & 240 \\ \hline 
			non-ISO      & 33 & 1 & 5 & 6 & 1 & 7 & 1 & 2 & 1 & 3 & 2 & 1 & 2 & 1 \\ \hline 
			co-iso.      & 2 & 0 & 0 & 0 & 0 & 0 & 0 & 0 & 1 & 0 & 0 & 0 & 0 & 1 \\ \hline 
			transitive      & 2 & 0 & 0 & 0 & 0 & 0 & 0 & 0 & 1 & 0 & 0 & 0 & 0 & 1 \\ \hline 
			st. co-iso.      & 2 & 0 & 0 & 0 & 0 & 0 & 0 & 0 & 1 & 0 & 0 & 0 & 0 & 1 \\ \hline 
			st. trans.      & 2 & 0 & 0 & 0 & 0 & 0 & 0 & 0 & 1 & 0 & 0 & 0 & 0 & 1 \\ \hline 
		\end{tabular}

	\caption{$5$-MOLR for $n=6$.}
	\label{6x5.}
\end{table}

\section{Sizes of autotopism groups of MOLR for $\mathbf{n=7}$}
\label{app_auto7}

\begin{table}[H]
	\small
		\setlength\tabcolsep{3pt}
		\begin{tabular}{|c|c||c|c|c|c|c|c|c|}
			\hline 
			2$\times$7 & $\Sigma$ & 1 & 2 & 4 & 14 & 24 & 28 & 48 \\ \hline 
			non-ISO      & 100 & 21 & 55 & 18 & 2 & 1 & 1 & 2 \\ \hline 
			co-iso.      & 42 & 3 & 16 & 18 & 2 & 0 & 1 & 2 \\ \hline 
			transitive      & 29 & 0 & 8 & 18 & 0 & 0 & 1 & 2 \\ \hline 
			st. co-iso.      & 42 & 3 & 16 & 18 & 2 & 0 & 1 & 2 \\ \hline 
			st. trans.      & 29 & 0 & 8 & 18 & 0 & 0 & 1 & 2 \\ \hline 
		\end{tabular}
		
		\begin{tabular}{|c|c||c|c|c|c|c|c|c|c|c|c|c|c|}
			\hline 
			3$\times$7 & $\Sigma$ & 1 & 2 & 3 & 4 & 6 & 7 & 14 & 21 & 24 & 28 & 42 & 72 \\ \hline 
			non-ISO      & \num{514162} & \num{508132} & \num{5880} & 48 & 65 & 23 & 4 & 4 & 2 & 1 & 1 & 1 & 1 \\ \hline 
			co-iso.      & \num{14464} & \num{10835} & \num{3524} & 6 & 65 & 23 & 1 & 4 & 2 & 1 & 1 & 1 & 1 \\ \hline 
			transitive      & \num{3549} & 0 & \num{3455} & 0 & 65 & 23 & 0 & 2 & 0 & 1 & 1 & 1 & 1 \\ \hline 
			st. co-iso.      & \num{7423} & \num{5017} & \num{2302} & 6 & 65 & 23 & 1 & 4 & 2 & 0 & 1 & 1 & 1 \\ \hline 
			st. trans.      & \num{2175} & 0 & \num{2082} & 0 & 65 & 23 & 0 & 2 & 0 & 0 & 1 & 1 & 1 \\ \hline 
		\end{tabular}
		
		\begin{tabular}{|c|c||c|c|c|c|c|c|c|c|c|c|c|}
			\hline 
			4$\times$7 & $\Sigma$ & 1 & 2 & 3 & 4 & 6 & 7 & 8 & 14 & 21 & 28 & 42 \\ \hline 
			non-ISO      & \num{49415812} & \num{49363791} & \num{51060} & 428 & 444 & 54 & 11 & 14 & 6 & 2 & 1 & 1 \\ \hline 
			co-iso.      & \num{65156} & \num{37639} & \num{27054} & 16 & 365 & 54 & 4 & 14 & 6 & 2 & 1 & 1 \\ \hline 
			transitive      & \num{27299} & 0 & \num{26867} & 0 & 361 & 54 & 0 & 14 & 1 & 0 & 1 & 1 \\ \hline 
			st. co-iso.      & \num{14960} & \num{4249} & \num{10418} & 16 & 205 & 54 & 0 & 11 & 3 & 2 & 1 & 1 \\ \hline 
			st. trans.      & \num{10029} & 0 & \num{9775} & 0 & 187 & 54 & 0 & 11 & 0 & 0 & 1 & 1 \\ \hline 
		\end{tabular}
		
		\begin{tabular}{|c|c||c|c|c|c|c|c|c|c|}
			\hline 
			5$\times$7 & $\Sigma$ & 1 & 2 & 4 & 5 & 7 & 10 & 14 & 28 \\ \hline 
			non-ISO      & \num{21290125} & \num{21243988} & \num{45872} & 227 & 10 & 6 & 12 & 9 & 1 \\ \hline 
			co-iso.      & \num{22432} & \num{3508} & \num{18672} & 227 & 1 & 6 & 8 & 9 & 1 \\ \hline 
			transitive      & \num{18836} & 0 & \num{18599} & 227 & 0 & 0 & 5 & 4 & 1 \\ \hline 
			st. co-iso.      & \num{4163} & 99 & \num{3935} & 121 & 0 & 0 & 5 & 2 & 1 \\ \hline 
			st. trans.      & \num{3923} & 0 & \num{3799} & 118 & 0 & 0 & 5 & 0 & 1 \\ \hline 
		\end{tabular}
		
		\begin{tabular}{|c|c||c|c|c|c|c|c|c|c|c|c|c|}
			\hline 
			6$\times$7 & $\Sigma$ & 1 & 2 & 3 & 4 & 5 & 6 & 7 & 10 & 12 & 42 & 84 \\ \hline 
			non-ISO      & \num{11582} & \num{10912} & 492 & 20 & 24 & 102 & 11 & 1 & 12 & 4 & 3 & 1 \\ \hline 
			co-iso.      & 409 & 2 & 345 & 0 & 24 & 8 & 9 & 1 & 12 & 4 & 3 & 1 \\ \hline 
			transitive      & 392 & 0 & 342 & 0 & 24 & 0 & 8 & 0 & 12 & 4 & 1 & 1 \\ \hline 
			st. co-iso.      & 91 & 0 & 61 & 0 & 4 & 2 & 7 & 0 & 12 & 2 & 2 & 1 \\ \hline 
			st. trans.      & 84 & 0 & 60 & 0 & 3 & 0 & 7 & 0 & 12 & 1 & 0 & 1 \\ \hline 
		\end{tabular}
		
		\begin{tabular}{|c|c||c|c|c|c|c|c|c|c|}
			\hline 
			7$\times$7 & $\Sigma$ & 2 & 3 & 6 & 12 & 21 & 42 & 294 & 588 \\ \hline 
			non-ISO      & 20 & 5 & 5 & 3 & 1 & 1 & 2 & 2 & 1 \\ \hline 
			co-iso.      & 9 & 1 & 0 & 2 & 1 & 0 & 2 & 2 & 1 \\ \hline 
			transitive      & 6 & 0 & 0 & 2 & 1 & 0 & 2 & 0 & 1 \\ \hline 
			st. co-iso.      & 6 & 0 & 0 & 1 & 1 & 0 & 1 & 2 & 1 \\ \hline 
			st. trans.      & 4 & 0 & 0 & 1 & 1 & 0 & 1 & 0 & 1 \\ \hline 
		\end{tabular}
		
	\caption{$2$-MOLR for $n=7$.}
	\label{7x2.}
\end{table}

\begin{table}[H]
\fontsize{9}{12}\selectfont
		\setlength\tabcolsep{3pt}
		\begin{tabular}{|c|c||c|c|c|c|c|c|c|c|}
			\hline 
			2$\times$7 & $\Sigma$ & 1 & 2 & 3 & 4 & 6 & 12 & 14 & 42 \\ \hline 
			non-ISO      & \num{2858} & \num{2300} & 512 & 3 & 28 & 9 & 2 & 3 & 1 \\ \hline 
			co-iso.      & 318 & 194 & 100 & 3 & 6 & 9 & 2 & 3 & 1 \\ \hline 
			transitive      & 15 & 0 & 0 & 3 & 0 & 9 & 2 & 0 & 1 \\ \hline 
			st. co-iso.      & 318 & 194 & 100 & 3 & 6 & 9 & 2 & 3 & 1 \\ \hline 
			st. trans.      & 15 & 0 & 0 & 3 & 0 & 9 & 2 & 0 & 1 \\ \hline 
		\end{tabular}
		
		\begin{tabular}{|c|c||c|c|c|c|c|c|c|c|c|c|c|}
			\hline 
			3$\times$7 & $\Sigma$ & 1 & 2 & 3 & 6 & 7 & 9 & 14 & 18 & 21 & 42 & 63 \\ \hline 
			non-ISO      & \num{65883453} & \num{65822447} & \num{60195} & 635 & 143 & 17 & 3 & 3 & 4 & 4 & 1 & 1 \\ \hline 
			co-iso.      & \num{49370} & \num{48126} & 566 & 542 & 116 & 4 & 3 & 3 & 4 & 4 & 1 & 1 \\ \hline 
			transitive      & 647 & 0 & 0 & 524 & 113 & 0 & 3 & 0 & 4 & 1 & 1 & 1 \\ \hline 
			st. co-iso.      & \num{13975} & \num{13397} & 305 & 189 & 64 & 4 & 3 & 3 & 4 & 4 & 1 & 1 \\ \hline 
			st. trans.      & 185 & 0 & 0 & 125 & 51 & 0 & 3 & 0 & 4 & 0 & 1 & 1 \\ \hline 
		\end{tabular}
		
		\begin{tabular}{|c|c||c|c|c|c|c|c|c|c|c|c|c|c|c|}
			\hline 
			4$\times$7 & $\Sigma$ & 1 & 2 & 3 & 4 & 6 & 7 & 9 & 12 & 14 & 18 & 21 & 42 & 63 \\ \hline 
			non-ISO      & \num{323112477} & \num{323002195} & \num{107997} & \num{1975} & 120 & 116 & 43 & 10 & 3 & 8 & 2 & 6 & 1 & 1 \\ \hline 
			co-iso.      & \num{2985} & \num{1232} & 147 & \num{1474} & 1 & 91 & 9 & 10 & 3 & 8 & 2 & 6 & 1 & 1 \\ \hline 
			transitive      & \num{1578} & 0 & 0 & \num{1468} & 0 & 90 & 0 & 10 & 3 & 0 & 2 & 3 & 1 & 1 \\ \hline 
			st. co-iso.      & 283 & 59 & 27 & 136 & 1 & 34 & 2 & 10 & 3 & 4 & 2 & 3 & 1 & 1 \\ \hline 
			st. trans.      & 160 & 0 & 0 & 113 & 0 & 30 & 0 & 10 & 3 & 0 & 2 & 0 & 1 & 1 \\ \hline 
		\end{tabular}
		
		\begin{tabular}{|c|c||c|c|c|c|c|c|c|c|c|}
			\hline 
			5$\times$7 & $\Sigma$ & 1 & 2 & 3 & 4 & 6 & 7 & 14 & 21 & 42 \\ \hline 
			non-ISO      & \num{55545} & \num{52981} & \num{2500} & 32 & 5 & 2 & 15 & 8 & 1 & 1 \\ \hline 
			co-iso.      & 111 & 31 & 21 & 32 & 0 & 2 & 15 & 8 & 1 & 1 \\ \hline 
			transitive      & 36 & 0 & 0 & 32 & 0 & 2 & 0 & 0 & 1 & 1 \\ \hline 
			st. co-iso.      & 8 & 0 & 0 & 2 & 0 & 2 & 0 & 3 & 0 & 1 \\ \hline 
			st. trans.      & 5 & 0 & 0 & 2 & 0 & 2 & 0 & 0 & 0 & 1 \\ \hline 
		\end{tabular}
		
		\begin{tabular}{|c|c||c|c|c|c|c|c|c|c|}
			\hline 
			6$\times$7 & $\Sigma$ & 1 & 2 & 3 & 4 & 6 & 12 & 42 & 126 \\ \hline 
			non-ISO      & 16 & 1 & 4 & 1 & 1 & 3 & 2 & 3 & 1 \\ \hline 
			co-iso.      & 11 & 0 & 1 & 0 & 1 & 3 & 2 & 3 & 1 \\ \hline 
			transitive      & 6 & 0 & 0 & 0 & 0 & 3 & 2 & 0 & 1 \\ \hline 
			st. co-iso.      & 4 & 0 & 0 & 0 & 0 & 0 & 0 & 3 & 1 \\ \hline 
			st. trans.      & 1 & 0 & 0 & 0 & 0 & 0 & 0 & 0 & 1 \\ \hline 
		\end{tabular}
		
		\begin{tabular}{|c|c||c|c|}
			\hline 
			7$\times$7 & $\Sigma$ & 294 & 882 \\ \hline 
			non-ISO      & 4 & 3 & 1 \\ \hline 
			co-iso.      & 4 & 3 & 1 \\ \hline 
			transitive      & 1 & 0 & 1 \\ \hline 
			st. co-iso.      & 4 & 3 & 1 \\ \hline 
			st. trans.      & 1 & 0 & 1 \\ \hline 
		\end{tabular}
		
	\caption{$3$-MOLR for $n=7$.}
	\label{7x3.}
\end{table}

\begin{table}[H]
		\setlength\tabcolsep{3pt}
		\begin{tabular}{|c|c||c|c|c|c|c|c|c|}
			\hline 
			2$\times$7 & $\Sigma$ & 1 & 2 & 4 & 8 & 14 & 16 & 28 \\ \hline 
			non-ISO      & \num{17609} & \num{15981} & \num{1545} & 64 & 12 & 2 & 4 & 1 \\ \hline 
			co-iso.      & 691 & 489 & 164 & 19 & 12 & 2 & 4 & 1 \\ \hline 
			transitive      & 21 & 0 & 0 & 5 & 12 & 0 & 4 & 0 \\ \hline 
			st. co-iso.      & 691 & 489 & 164 & 19 & 12 & 2 & 4 & 1 \\ \hline 
			st. trans.      & 21 & 0 & 0 & 5 & 12 & 0 & 4 & 0 \\ \hline 
		\end{tabular}
		
		\begin{tabular}{|c|c||c|c|c|c|c|c|c|c|c|c|c|}
			\hline 
			3$\times$7 & $\Sigma$ & 1 & 2 & 3 & 4 & 6 & 7 & 12 & 14 & 21 & 28 & 42 \\ \hline 
			non-ISO      & \num{35469948} & \num{35420362} & \num{48685} & 626 & 226 & 21 & 12 & 1 & 7 & 6 & 1 & 1 \\ \hline 
			co-iso.      & \num{1622} & 949 & 539 & 9 & 110 & 4 & 1 & 1 & 4 & 3 & 1 & 1 \\ \hline 
			transitive      & 110 & 0 & 0 & 0 & 109 & 0 & 0 & 1 & 0 & 0 & 0 & 0 \\ \hline 
			st. co-iso.      & 585 & 251 & 243 & 9 & 67 & 4 & 1 & 1 & 4 & 3 & 1 & 1 \\ \hline 
			st. trans.      & 48 & 0 & 0 & 0 & 47 & 0 & 0 & 1 & 0 & 0 & 0 & 0 \\ \hline 
		\end{tabular}
		
		\begin{tabular}{|c|c||c|c|c|c|c|c|c|c|c|c|c|c|c|c|}
			\hline 
			4$\times$7 & $\Sigma$ & 1 & 2 & 3 & 4 & 6 & 7 & 8 & 14 & 16 & 21 & 28 & 42 & 56 & 84 \\ \hline 
			non-ISO      & \num{68659} & \num{67073} & \num{1354} & 74 & 110 & 2 & 20 & 12 & 5 & 1 & 3 & 2 & 1 & 1 & 1 \\ \hline 
			co-iso.      & 84 & 3 & 5 & 0 & 52 & 0 & 1 & 12 & 3 & 1 & 2 & 2 & 1 & 1 & 1 \\ \hline 
			transitive      & 67 & 0 & 0 & 0 & 51 & 0 & 0 & 12 & 0 & 1 & 0 & 1 & 0 & 1 & 1 \\ \hline 
			st. co-iso.      & 12 & 0 & 1 & 0 & 2 & 0 & 0 & 2 & 2 & 1 & 2 & 1 & 1 & 0 & 0 \\ \hline 
			st. trans.      & 1 & 0 & 0 & 0 & 0 & 0 & 0 & 1 & 0 & 0 & 0 & 0 & 0 & 0 & 0 \\ \hline 
		\end{tabular}
		
		\begin{tabular}{|c|c||c|c|c|c|c|c|c|c|}
			\hline 
			5$\times$7 & $\Sigma$ & 1 & 2 & 4 & 7 & 8 & 14 & 28 & 56 \\ \hline 
			non-ISO      & 204 & 96 & 41 & 32 & 6 & 18 & 6 & 4 & 1 \\ \hline 
			co-iso.      & 67 & 1 & 0 & 31 & 6 & 18 & 6 & 4 & 1 \\ \hline 
			transitive      & 53 & 0 & 0 & 31 & 0 & 18 & 0 & 3 & 1 \\ \hline 
			st. co-iso.      & 3 & 0 & 0 & 0 & 0 & 0 & 2 & 1 & 0 \\ \hline 
		\end{tabular}
		
		\begin{tabular}{|c|c||c|c|c|c|c|}
			\hline 
			6$\times$7 & $\Sigma$ & 4 & 6 & 24 & 42 & 84 \\ \hline 
			non-ISO      & 7 & 1 & 1 & 2 & 2 & 1 \\ \hline 
			co-iso.      & 7 & 1 & 1 & 2 & 2 & 1 \\ \hline 
			transitive      & 3 & 1 & 0 & 2 & 0 & 0 \\ \hline 
			st. co-iso.      & 3 & 0 & 0 & 0 & 2 & 1 \\ \hline 
		\end{tabular}
		
		\begin{tabular}{|c|c||c|c|}
			\hline 
			7$\times$7 & $\Sigma$ & 294 & 588 \\ \hline 
			non-ISO      & 3 & 2 & 1 \\ \hline 
			co-iso.      & 3 & 2 & 1 \\ \hline 
			st. co-iso.      & 3 & 2 & 1 \\ \hline 
		\end{tabular}
		
	\caption{$4$-MOLR for $n=7$.}
	\label{7x4.}
\end{table}

\begin{table}[H]
		
		\begin{tabular}{|c|c||c|c|c|c|c|c|c|}
			\hline 
			2$\times$7 & $\Sigma$ & 1 & 2 & 4 & 5 & 8 & 10 & 14 \\ \hline 
			non-ISO      & \num{10626} & \num{9590} & 957 & 48 & 4 & 24 & 2 & 1 \\ \hline 
			co-iso.      & 176 & 122 & 37 & 4 & 4 & 6 & 2 & 1 \\ \hline 
			transitive      & 6 & 0 & 0 & 0 & 4 & 0 & 2 & 0 \\ \hline 
			st. co-iso.      & 176 & 122 & 37 & 4 & 4 & 6 & 2 & 1 \\ \hline 
			st. trans.      & 6 & 0 & 0 & 0 & 4 & 0 & 2 & 0 \\ \hline 
		\end{tabular}
		
		\begin{tabular}{|c|c||c|c|c|c|c|c|c|c|c|}
			\hline 
			3$\times$7 & $\Sigma$ & 1 & 2 & 3 & 5 & 6 & 7 & 10 & 14 & 21 \\ \hline 
			non-ISO      & \num{22982} & \num{21848} & \num{1039} & 39 & 30 & 7 & 1 & 12 & 1 & 5 \\ \hline 
			co-iso.      & 49 & 3 & 1 & 0 & 30 & 0 & 0 & 12 & 1 & 2 \\ \hline 
			transitive      & 42 & 0 & 0 & 0 & 30 & 0 & 0 & 12 & 0 & 0 \\ \hline 
			st. co-iso.      & 48 & 2 & 1 & 0 & 30 & 0 & 0 & 12 & 1 & 2 \\ \hline 
			st. trans.      & 42 & 0 & 0 & 0 & 30 & 0 & 0 & 12 & 0 & 0 \\ \hline 
		\end{tabular}
		
		\begin{tabular}{|c|c||c|c|c|c|c|c|}
			\hline 
			4$\times$7 & $\Sigma$ & 1 & 2 & 4 & 8 & 14 & 21 \\ \hline 
			non-ISO      & 19 & 2 & 4 & 8 & 3 & 1 & 1 \\ \hline 
			co-iso.      & 2 & 0 & 0 & 0 & 0 & 1 & 1 \\ \hline 
			st. co-iso.      & 2 & 0 & 0 & 0 & 0 & 1 & 1 \\ \hline 
		\end{tabular}
		
		\begin{tabular}{|c|c||c|c|c|}
			\hline 
			5$\times$7 & $\Sigma$ & 4 & 14 & 20 \\ \hline 
			non-ISO      & 5 & 2 & 1 & 2 \\ \hline 
			co-iso.      & 4 & 1 & 1 & 2 \\ \hline 
			transitive      & 2 & 0 & 0 & 2 \\ \hline 
			st. co-iso.      & 1 & 0 & 1 & 0 \\ \hline 
		\end{tabular}
		
		\begin{tabular}{|c|c||c|c|c|c|}
			\hline 
			6$\times$7 & $\Sigma$ & 20 & 24 & 42 & 120 \\ \hline 
			non-ISO      & 5 & 1 & 1 & 1 & 2 \\ \hline 
			co-iso.      & 5 & 1 & 1 & 1 & 2 \\ \hline 
			transitive      & 3 & 1 & 0 & 0 & 2 \\ \hline 
			st. co-iso.      & 1 & 0 & 0 & 1 & 0 \\ \hline 
		\end{tabular}
		
		\begin{tabular}{|c|c||c|}
			\hline 
			7$\times$7 & $\Sigma$ & 294 \\ \hline 
			non-ISO      & 1 & 1 \\ \hline 
			co-iso.      & 1 & 1 \\ \hline 
			st. co-iso.      & 1 & 1 \\ \hline 
		\end{tabular}
		
	\caption{$5$-MOLR for $n=7$.}
	\label{7x5.}
\end{table}

\begin{table}[H]
		\setlength\tabcolsep{3pt}
		\begin{tabular}{|c|c||c|c|c|c|c|c|c|c|c|c|c|c|}
			\hline 
			2$\times$7 & $\Sigma$ & 1 & 2 & 3 & 4 & 5 & 6 & 8 & 10 & 12 & 16 & 48 & 84 \\ \hline 
			non-ISO      & \num{1895} & \num{1505} & 328 & 2 & 29 & 4 & 5 & 14 & 2 & 1 & 2 & 2 & 1 \\ \hline 
			co-iso.      & 26 & 7 & 7 & 1 & 1 & 1 & 2 & 0 & 1 & 1 & 2 & 2 & 1 \\ \hline 
			transitive      & 5 & 0 & 0 & 0 & 0 & 0 & 1 & 0 & 0 & 1 & 0 & 2 & 1 \\ \hline 
			st. co-iso.      & 26 & 7 & 7 & 1 & 1 & 1 & 2 & 0 & 1 & 1 & 2 & 2 & 1 \\ \hline 
			st. trans.      & 5 & 0 & 0 & 0 & 0 & 0 & 1 & 0 & 0 & 1 & 0 & 2 & 1 \\ \hline 
		\end{tabular}
		
		\begin{tabular}{|c|c||c|c|c|c|c|c|c|c|}
			\hline 
			3$\times$7 & $\Sigma$ & 1 & 2 & 3 & 6 & 12 & 21 & 84 & 126 \\ \hline 
			non-ISO      & 23 & 8 & 4 & 1 & 5 & 1 & 1 & 1 & 2 \\ \hline 
			co-iso.      & 7 & 0 & 0 & 0 & 3 & 1 & 0 & 1 & 2 \\ \hline 
			transitive      & 7 & 0 & 0 & 0 & 3 & 1 & 0 & 1 & 2 \\ \hline 
			st. co-iso.      & 6 & 0 & 0 & 0 & 2 & 1 & 0 & 1 & 2 \\ \hline 
			st. trans.      & 4 & 0 & 0 & 0 & 0 & 1 & 0 & 1 & 2 \\ \hline 
		\end{tabular}
		
		\begin{tabular}{|c|c||c|c|}
			\hline 
			4$\times$7 & $\Sigma$ & 84 & 126 \\ \hline 
			non-ISO      & 2 & 1 & 1 \\ \hline 
			co-iso.      & 2 & 1 & 1 \\ \hline 
			transitive      & 2 & 1 & 1 \\ \hline 
			st. co-iso.      & 2 & 1 & 1 \\ \hline 
			st. trans.      & 2 & 1 & 1 \\ \hline 
		\end{tabular}
		
		\begin{tabular}{|c|c||c|}
			\hline 
			5$\times$7 & $\Sigma$ & 84 \\ \hline 
			non-ISO      & 1 & 1 \\ \hline 
			co-iso.      & 1 & 1 \\ \hline 
			transitive      & 1 & 1 \\ \hline 
			st. co-iso.      & 1 & 1 \\ \hline 
			st. trans.      & 1 & 1 \\ \hline 
		\end{tabular}
		
		\begin{tabular}{|c|c||c|}
			\hline 
			6$\times$7 & $\Sigma$ & 252 \\ \hline 
			non-ISO      & 1 & 1 \\ \hline 
			co-iso.      & 1 & 1 \\ \hline 
			transitive      & 1 & 1 \\ \hline 
			st. co-iso.      & 1 & 1 \\ \hline 
			st. trans.      & 1 & 1 \\ \hline 
		\end{tabular}
		
		\begin{tabular}{|c|c||c|}
			\hline 
			7$\times$7 & $\Sigma$ & \num{1764} \\ \hline 
			non-ISO      & 1 & 1 \\ \hline 
			co-iso.      & 1 & 1 \\ \hline 
			transitive      & 1 & 1 \\ \hline 
			st. co-iso.      & 1 & 1 \\ \hline 
			st. trans.      & 1 & 1 \\ \hline 
		\end{tabular}
		
	\caption{$6$-MOLR for $n=7$.}
	\label{7x6.}
\end{table}

\section{Sizes of autotopism groups of MOLR for $\mathbf{n=8}$}
\label{app_auto8}

\begin{landscape}
\begin{table}[H]
		\setlength\tabcolsep{5pt}
		\begin{tabular}{|c|c||c|c|c|c|c|c|c|c|c|c|c|}
			\hline 
			2$\times$8 & $\Sigma$ & 1 & 2 & 4 & 8 & 12 & 16 & 30 & 32 & 60 & 64 & 128 \\ \hline 
			st. co-iso.      & 186 & 52 & 65 & 43 & 10 & 2 & 2 & 1 & 6 & 1 & 2 & 2 \\ \hline 
			st. trans.      & 99 & 0 & 35 & 39 & 10 & 2 & 2 & 0 & 6 & 1 & 2 & 2 \\ \hline 
		\end{tabular}
		
		\begin{tabular}{|c|c||c|c|c|c|c|c|c|c|c|c|c|c|c|c|}
			\hline 
			3$\times$8 & $\Sigma$ & 1 & 2 & 3 & 4 & 6 & 8 & 12 & 16 & 30 & 32 & 48 & 60 & 64 & 192 \\ \hline 
			st. co-iso.      & \num{446443} & \num{394387} & \num{50556} & 11 & \num{1311} & 28 & 101 & 22 & 12 & 2 & 8 & 1 & 1 & 1 & 2 \\ \hline 
			st. trans.      & \num{45429} & 0 & \num{43990} & 0 & \num{1265} & 27 & 99 & 22 & 12 & 1 & 8 & 1 & 1 & 1 & 2 \\ \hline 
		\end{tabular}
		
		\begin{tabular}{|c|c||c|c|c|c|c|c|c|c|c|c|c|c|c|c|}
			\hline 
			4$\times$8 & $\Sigma$ & 1 & 2 & 3 & 4 & 6 & 8 & 12 & 16 & 24 & 32 & 48 & 64 & 192 & 768 \\ \hline 
			st. co-iso.      & \num{4432284} & \num{3196674} & \num{1222057} & 168 & \num{11703} & 343 & \num{1002} & 50 & 194 & 11 & 54 & 11 & 14 & 1 & 2 \\ \hline 
			st. trans.      & \num{1097655} & 0 & \num{1086038} & 0 & \num{10123} & 317 & 857 & 50 & 181 & 9 & 52 & 11 & 14 & 1 & 2 \\ \hline 
		\end{tabular}
		
		\begin{tabular}{|c|c||c|c|c|c|c|c|}
			\hline 
			5$\times$8 & $\Sigma$ & 1 & 2 & 4 & 8 & 16 & 32 \\ \hline 
			st. co-iso.      & \num{3826527} & \num{884803} & \num{2929324} & \num{11629} & 702 & 65 & 4 \\ \hline 
			st. trans.      & \num{2569679} & 0 & \num{2558291} & \num{10693} & 629 & 62 & 4 \\ \hline 
		\end{tabular}
		
		\begin{tabular}{|c|c||c|c|c|c|c|c|c|c|c|c|}
			\hline 
			6$\times$8 & $\Sigma$ & 1 & 2 & 3 & 4 & 6 & 8 & 12 & 16 & 32 & 96 \\ \hline 
			st. co-iso.      & \num{242732} & \num{7954} & \num{227619} & 2 & \num{6289} & 46 & 611 & 23 & 181 & 6 & 1 \\ \hline 
			st. trans.      & \num{206612} & 0 & \num{200153} & 0 & \num{5745} & 35 & 503 & 22 & 147 & 6 & 1 \\ \hline 
		\end{tabular}
		
		\begin{tabular}{|c|c||c|c|c|c|c|c|c|c|c|c|c|c|c|}
			\hline 
			7$\times$8 & $\Sigma$ & 1 & 2 & 3 & 4 & 6 & 8 & 12 & 16 & 24 & 32 & 48 & 56 & 96 \\ \hline 
			st. co-iso.      & 484 & 28 & 175 & 5 & 109 & 72 & 48 & 29 & 9 & 2 & 2 & 2 & 1 & 2 \\ \hline 
			st. trans.      & 305 & 0 & 92 & 0 & 82 & 64 & 24 & 29 & 7 & 2 & 2 & 1 & 0 & 2 \\ \hline 
		\end{tabular}
		
		\begin{tabular}{|c|c||c|c|c|c|c|c|}
			\hline 
			8$\times$8 & $\Sigma$ & 4 & 8 & 16 & 32 & 64 & 448 \\ \hline 
			st. co-iso.      & 70 & 9 & 19 & 22 & 8 & 11 & 1 \\ \hline 
			st. trans.      & 13 & 0 & 2 & 4 & 4 & 3 & 0 \\ \hline 
		\end{tabular}
		
	\caption{$2$-MOLR for $n=8$.}
	\label{8x2.}
\end{table}

\begin{table}[H]
		\setlength\tabcolsep{3pt}
		\begin{tabular}{|c|c||c|c|c|c|c|c|c|c|c|c|c|c|c|}
			\hline 
			2$\times$8 & $\Sigma$ & 1 & 2 & 3 & 4 & 6 & 8 & 12 & 16 & 24 & 32 & 48 & 96 & 384 \\ \hline 
			co-iso.      & \num{11565} & \num{10583} & 803 & 24 & 83 & 22 & 20 & 14 & 6 & 2 & 4 & 2 & 1 & 1 \\ \hline 
			transitive      & 66 & 0 & 0 & 24 & 0 & 22 & 0 & 14 & 0 & 2 & 0 & 2 & 1 & 1 \\ \hline 
			st. co-iso.      & \num{11565} & \num{10583} & 803 & 24 & 83 & 22 & 20 & 14 & 6 & 2 & 4 & 2 & 1 & 1 \\ \hline 
			st. trans.      & 66 & 0 & 0 & 24 & 0 & 22 & 0 & 14 & 0 & 2 & 0 & 2 & 1 & 1 \\ \hline 
		\end{tabular}
		
		\begin{tabular}{|c|c||c|c|c|c|c|c|c|c|c|c|c|c|c|c|}
			\hline 
			3$\times$8 & $\Sigma$ & 1 & 2 & 3 & 4 & 6 & 8 & 12 & 16 & 18 & 24 & 36 & 48 & 144 & 576 \\ \hline 
			st. co-iso.      & \num{9144025} & \num{9121524} & \num{13878} & \num{7463} & 479 & 564 & 27 & 64 & 8 & 4 & 3 & 3 & 5 & 1 & 2 \\ \hline 
			st. trans.      & \num{7627} & 0 & 0 & \num{7067} & 0 & 485 & 0 & 57 & 0 & 4 & 3 & 3 & 5 & 1 & 2 \\ \hline 
		\end{tabular}
		
		\begin{tabular}{|c|c||c|c|c|c|c|c|c|c|c|c|c|c|c|c|c|c|c|c|c|c|}
			\hline 
			4$\times$8 & $\Sigma$ & 1 & 2 & 3 & 4 & 6 & 8 & 9 & 12 & 16 & 18 & 24 & 32 & 36 & 48 & 64 & 72 & 96 & 144 & 192 & \num{2304} \\ \hline 
			st. co-iso.      & \num{178502} & \num{91562} & \num{40127} & \num{41876} & \num{2464} & \num{1665} & 397 & 38 & 177 & 94 & 26 & 39 & 8 & 3 & 10 & 5 & 2 & 3 & 1 & 3 & 2 \\ \hline 
			st. trans.      & \num{41505} & 0 & 0 & \num{39650} & 0 & \num{1566} & 0 & 38 & 168 & 0 & 26 & 36 & 0 & 3 & 7 & 0 & 2 & 3 & 1 & 3 & 2 \\ \hline 
		\end{tabular}
		
		\begin{tabular}{|c|c||c|c|c|c|c|c|c|c|c|c|}
			\hline 
			5$\times$8 & $\Sigma$ & 1 & 2 & 3 & 4 & 6 & 8 & 12 & 16 & 24 & 48 \\ \hline 
			st. co-iso.      & 628 & 34 & 142 & 311 & 59 & 34 & 26 & 5 & 5 & 10 & 2 \\ \hline 
			st. trans.      & 75 & 0 & 0 & 38 & 0 & 25 & 0 & 3 & 0 & 7 & 2 \\ \hline 
		\end{tabular}
		
		\begin{tabular}{|c|c||c|c|c|c|c|c|c|c|c|}
			\hline 
			6$\times$8 & $\Sigma$ & 2 & 3 & 4 & 6 & 8 & 12 & 16 & 24 & 48 \\ \hline 
			st. co-iso.      & 111 & 18 & 8 & 42 & 8 & 8 & 12 & 8 & 2 & 5 \\ \hline 
			st. trans.      & 32 & 0 & 8 & 0 & 8 & 0 & 12 & 0 & 2 & 2 \\ \hline 
		\end{tabular}
		
		\begin{tabular}{|c|c||c|c|c|c|c|}
			\hline 
			7$\times$8 & $\Sigma$ & 8 & 9 & 24 & 56 & 168 \\ \hline 
			st. co-iso.      & 10 & 2 & 3 & 2 & 1 & 2 \\ \hline 
			st. trans.      & 6 & 0 & 3 & 1 & 0 & 2 \\ \hline 
		\end{tabular}
		
		\begin{tabular}{|c|c||c|c|c|c|c|}
			\hline 
			8$\times$8 & $\Sigma$ & 48 & 64 & 192 & 448 & \num{1344} \\ \hline 
			st. co-iso.      & 7 & 1 & 2 & 1 & 1 & 2 \\ \hline 
			st. trans.      & 3 & 1 & 0 & 0 & 0 & 2 \\ \hline 
		\end{tabular}
		
	\caption{$3$-MOLR for $n=8$.}
	\label{8x3.}
\end{table}

\begin{table}[H]
		
		\begin{tabular}{|c|c||c|c|c|c|c|c|c|c|c|c|c|c|c|c|c|}
			\hline 
			2$\times$8 & $\Sigma$ & 1 & 2 & 3 & 4 & 6 & 8 & 12 & 16 & 24 & 32 & 48 & 64 & 96 & 128 & 384 \\ \hline 
			st. co-iso.      & \num{216950} & \num{212259} & \num{4241} & 20 & 313 & 8 & 52 & 2 & 31 & 1 & 9 & 6 & 3 & 1 & 3 & 1 \\ \hline 
			st. trans.      & 152 & 0 & 0 & 0 & 98 & 0 & 17 & 0 & 20 & 1 & 4 & 5 & 3 & 0 & 3 & 1 \\ \hline 
		\end{tabular}
		
		\begin{tabular}{|c|c||c|c|c|c|c|c|c|c|c|c|c|}
			\hline 
			3$\times$8 & $\Sigma$ & 1 & 2 & 3 & 4 & 6 & 8 & 12 & 16 & 24 & 32 & 48 \\ \hline 
			st. co-iso.      & \num{1648723} & \num{1596362} & \num{45732} & 46 & \num{6203} & 20 & 279 & 5 & 58 & 6 & 8 & 4 \\ \hline 
			st. trans.      & \num{4284} & 0 & 0 & 0 & \num{4028} & 0 & 196 & 2 & 46 & 2 & 8 & 2 \\ \hline 
		\end{tabular}
		
		\begin{tabular}{|c|c||c|c|c|c|c|c|c|c|c|c|c|c|c|}
			\hline 
			4$\times$8 & $\Sigma$ & 1 & 2 & 3 & 4 & 6 & 8 & 12 & 16 & 24 & 32 & 48 & 64 & 96 \\ \hline 
			st. co-iso.      & \num{3547} & 296 & \num{1420} & 14 & \num{1273} & 5 & 340 & 5 & 106 & 12 & 45 & 1 & 26 & 4 \\ \hline 
			st. trans.      & 712 & 0 & 0 & 0 & 510 & 0 & 124 & 4 & 32 & 6 & 15 & 0 & 21 & 0 \\ \hline 
		\end{tabular}
		
		\begin{tabular}{|c|c||c|c|c|c|c|c|c|c|}
			\hline 
			5$\times$8 & $\Sigma$ & 2 & 4 & 8 & 16 & 24 & 32 & 48 & 64 \\ \hline 
			st. co-iso.      & 58 & 2 & 10 & 25 & 9 & 4 & 4 & 2 & 2 \\ \hline 
			st. trans.      & 20 & 0 & 3 & 5 & 3 & 1 & 4 & 2 & 2 \\ \hline 
		\end{tabular}
		
		\begin{tabular}{|c|c||c|c|c|}
			\hline 
			6$\times$8 & $\Sigma$ & 16 & 24 & 48 \\ \hline 
			st. co-iso.      & 4 & 1 & 1 & 2 \\ \hline 
		\end{tabular}
		
		\begin{tabular}{|c|c||c|c|}
			\hline 
			7$\times$8 & $\Sigma$ & 56 & 168 \\ \hline 
			st. co-iso.      & 3 & 1 & 2 \\ \hline 
		\end{tabular}
		
		\begin{tabular}{|c|c||c|c|}
			\hline 
			8$\times$8 & $\Sigma$ & 448 & \num{1344} \\ \hline 
			st. co-iso.      & 3 & 1 & 2 \\ \hline 
		\end{tabular}
		
	\caption{$4$-MOLR for $n=8$.}
	\label{8x4.}
\end{table}
\begin{table}[H]
		
		\begin{tabular}{|c|c||c|c|c|c|c|c|c|c|c|c|c|}
			\hline 
			2$\times$8 & $\Sigma$ & 1 & 2 & 4 & 5 & 8 & 10 & 16 & 20 & 32 & 64 & 128 \\ \hline 
			st. co-iso.      & \num{509622} & \num{505439} & \num{3896} & 185 & 9 & 52 & 8 & 21 & 2 & 5 & 3 & 2 \\ \hline 
			st. trans.      & 19 & 0 & 0 & 0 & 9 & 0 & 8 & 0 & 2 & 0 & 0 & 0 \\ \hline 
		\end{tabular}
		
		\begin{tabular}{|c|c||c|c|c|c|c|c|c|c|c|c|}
			\hline 
			3$\times$8 & $\Sigma$ & 1 & 2 & 3 & 4 & 6 & 8 & 12 & 16 & 24 & 48 \\ \hline 
			st. co-iso.      & \num{2652} & \num{2112} & 419 & 2 & 75 & 13 & 15 & 2 & 9 & 2 & 3 \\ \hline 
		\end{tabular}
		
		\begin{tabular}{|c|c||c|c|c|c|c|c|c|c|c|c|c|c|c|c|}
			\hline 
			4$\times$8 & $\Sigma$ & 1 & 2 & 3 & 4 & 6 & 8 & 12 & 16 & 24 & 32 & 48 & 64 & 96 & 192 \\ \hline 
			st. co-iso.      & 267 & 129 & 30 & 12 & 30 & 1 & 23 & 3 & 17 & 2 & 11 & 1 & 4 & 2 & 2 \\ \hline 
		\end{tabular}
		
		\begin{tabular}{|c|c||c|c|}
			\hline 
			5$\times$8 & $\Sigma$ & 8 & 16 \\ \hline 
			st. co-iso.      & 2 & 1 & 1 \\ \hline 
		\end{tabular}
		
		\begin{tabular}{|c|c||c|}
			\hline 
			6$\times$8 & $\Sigma$ & 16 \\ \hline 
			st. co-iso.      & 1 & 1 \\ \hline 
		\end{tabular}
		
		\begin{tabular}{|c|c||c|}
			\hline 
			7$\times$8 & $\Sigma$ & 56 \\ \hline 
			st. co-iso.      & 1 & 1 \\ \hline 
		\end{tabular}
		
		\begin{tabular}{|c|c||c|}
			\hline 
			8$\times$8 & $\Sigma$ & 448 \\ \hline 
			st. co-iso.      & 1 & 1 \\ \hline 
		\end{tabular}
		
	\caption{$5$-MOLR for $n=8$.}
	\label{8x5.}
\end{table}

\begin{table}[H]
		\setlength\tabcolsep{5pt}
		\begin{tabular}{|c|c||c|c|c|c|c|c|c|c|c|c|c|c|c|c|c|c|c|}
			\hline 
			2$\times$8 & $\Sigma$ & 1 & 2 & 3 & 4 & 5 & 6 & 8 & 10 & 12 & 16 & 24 & 32 & 48 & 64 & 96 & 128 & 384 \\ \hline 
			st. co-iso.      & \num{91013} & \num{89017} & \num{1597} & 59 & 138 & 4 & 71 & 35 & 8 & 37 & 22 & 3 & 7 & 7 & 2 & 2 & 2 & 2 \\ \hline 
			st. trans.      & 109 & 0 & 0 & 0 & 0 & 0 & 60 & 0 & 0 & 36 & 0 & 3 & 0 & 6 & 0 & 2 & 0 & 2 \\ \hline 
		\end{tabular}
		
		\begin{tabular}{|c|c||c|c|c|c|c|c|c|c|c|c|c|c|c|}
			\hline 
			3$\times$8 & $\Sigma$ & 1 & 2 & 3 & 4 & 6 & 8 & 12 & 18 & 24 & 36 & 48 & 96 & 288 \\ \hline 
			st. co-iso.      & 975 & 3 & 10 & 18 & 7 & 797 & 3 & 120 & 2 & 5 & 5 & 1 & 3 & 1 \\ \hline 
			st. trans.      & 908 & 0 & 0 & 0 & 0 & 775 & 0 & 119 & 0 & 4 & 5 & 1 & 3 & 1 \\ \hline 
		\end{tabular}
		
		\begin{tabular}{|c|c||c|c|c|c|c|c|c|c|c|c|c|}
			\hline 
			4$\times$8 & $\Sigma$ & 3 & 6 & 9 & 12 & 18 & 24 & 32 & 36 & 48 & 96 & 144 \\ \hline 
			st. co-iso.      & 155 & 2 & 122 & 2 & 14 & 8 & 1 & 1 & 1 & 1 & 2 & 1 \\ \hline 
			st. trans.      & 146 & 0 & 122 & 0 & 13 & 8 & 0 & 0 & 1 & 1 & 0 & 1 \\ \hline 
		\end{tabular}
		
		\begin{tabular}{|c|c||c|}
			\hline 
			5$\times$8 & $\Sigma$ & 24 \\ \hline 
			st. co-iso.      & 1 & 1 \\ \hline 
		\end{tabular}
		
		\begin{tabular}{|c|c||c|}
			\hline 
			6$\times$8 & $\Sigma$ & 48 \\ \hline 
			st. co-iso.      & 1 & 1 \\ \hline 
		\end{tabular}
		
		\begin{tabular}{|c|c||c|}
			\hline 
			7$\times$8 & $\Sigma$ & 168 \\ \hline 
			st. co-iso.      & 1 & 1 \\ \hline 
		\end{tabular}
		
		\begin{tabular}{|c|c||c|}
			\hline 
			8$\times$8 & $\Sigma$ & \num{1344} \\ \hline 
			st. co-iso.      & 1 & 1 \\ \hline 
		\end{tabular}
		
	\caption{$6$-MOLR for $n=8$.}
	\label{8x6.}
\end{table}

\begin{table}[H]
		
		\begin{tabular}{|c|c||c|c|c|c|c|c|c|c|c|c|c|c|c|c|c|c|}
			\hline 
			2$\times$8 & $\Sigma$ & 1 & 2 & 3 & 4 & 5 & 6 & 7 & 8 & 10 & 32 & 42 & 48 & 84 & 128 & 192 & \num{2688} \\ \hline 
			st. co-iso.      & \num{4538} & \num{4439} & 67 & 9 & 3 & 2 & 3 & 1 & 2 & 4 & 1 & 2 & 1 & 1 & 1 & 1 & 1 \\ \hline 
			st. trans.      & 5 & 0 & 0 & 0 & 0 & 0 & 0 & 1 & 0 & 0 & 0 & 2 & 0 & 1 & 0 & 0 & 1 \\ \hline 
		\end{tabular}
		
		\begin{tabular}{|c|c||c|c|}
			\hline 
			3$\times$8 & $\Sigma$ & 21 & 168 \\ \hline 
			st. co-iso.      & 2 & 1 & 1 \\ \hline 
			st. trans.      & 2 & 1 & 1 \\ \hline 
		\end{tabular}
		
		\begin{tabular}{|c|c||c|c|}
			\hline 
			4$\times$8 & $\Sigma$ & 168 & 672 \\ \hline 
			st. co-iso.      & 2 & 1 & 1 \\ \hline 
			st. trans.      & 2 & 1 & 1 \\ \hline 
		\end{tabular}
		
		\begin{tabular}{|c|c||c|}
			\hline 
			5$\times$8 & $\Sigma$ & 168 \\ \hline 
			st. co-iso.      & 1 & 1 \\ \hline 
			st. trans.      & 1 & 1 \\ \hline 
		\end{tabular}
		
		\begin{tabular}{|c|c||c|}
			\hline 
			6$\times$8 & $\Sigma$ & 336 \\ \hline 
			st. co-iso.      & 1 & 1 \\ \hline 
			st. trans.      & 1 & 1 \\ \hline 
		\end{tabular}
		
		\begin{tabular}{|c|c||c|}
			\hline 
			7$\times$8 & $\Sigma$ & \num{1176} \\ \hline 
			st. co-iso.      & 1 & 1 \\ \hline 
			st. trans.      & 1 & 1 \\ \hline 
		\end{tabular}
		
		\begin{tabular}{|c|c||c|}
			\hline 
			8$\times$8 & $\Sigma$ & \num{9408} \\ \hline 
			st. co-iso.      & 1 & 1 \\ \hline 
			st. trans.      & 1 & 1 \\ \hline 
		\end{tabular}
		
	\caption{$7$-MOLR for $n=8$.}
	\label{8x7.}
\end{table}
\end{landscape}

\section{Sizes of autotopism groups of MOLR for $\mathbf{n=9}$}
\label{app_auto9}

\begin{landscape}

\begin{table}[H]

		\begin{tabular}{|c|c||c|c|c|c|c|c|c|c|c|}
			\hline 
			2$\times$9 & $\Sigma$ & 4 & 6 & 12 & 16 & 24 & 36 & 72 & 80 & 648 \\ \hline 
			st. trans.      & 126 & 95 & 3 & 13 & 4 & 3 & 3 & 2 & 2 & 1 \\ \hline 
		\end{tabular}
		
		\begin{tabular}{|c|c||c|c|c|c|c|c|c|c|c|c|c|c|c|}
			\hline 
			3$\times$9 & $\Sigma$ & 2 & 4 & 6 & 8 & 12 & 18 & 24 & 36 & 40 & 54 & 72 & 108 & \num{1944} \\ \hline 
			st. trans.      & \num{1418577} & \num{1413987} & \num{3905} & 511 & 8 & 115 & 19 & 8 & 13 & 3 & 2 & 3 & 2 & 1 \\ \hline 
		\end{tabular}
		\setlength\tabcolsep{4pt}			
		\begin{tabular}{|c|c||c|c|c|c|c|c|c|c|c|c|c|c|c|c|}
			\hline 
			4$\times$9 & $\Sigma$ & 2 & 4 & 6 & 8 & 12 & 16 & 18 & 24 & 36 & 54 & 72 & 80 & 108 & 144 \\ \hline 
			st. trans.      & \num{560524587} & \num{560437428} & \num{82404} & \num{3520} & 958 & 220 & 15 & 19 & 12 & 5 & 2 & 1 & 1 & 1 & 1 \\ \hline 
		\end{tabular}
	
		\begin{tabular}{|c|c||c|c|c|c|c|c|c|c|c|c|c|c|c|}
			\hline 
			5$\times$9 & $\Sigma$ & 2 & 4 & 6 & 8 & 12 & 16 & 18 & 24 & 36 & 54 & 72 & 108 & 144 \\ \hline 
			st. trans.  & 20019499500 & 20019088250 & \num{400369} & \num{9228} & \num{1223} & 370 & 6 & 29 & 9 & 11 & 2 & 1 & 1 & 1 \\ \hline 
		\end{tabular}
		
		\begin{tabular}{|c|c||c|c|c|c|c|c|c|c|c|c|c|c|c|}
			\hline 
			6$\times$9 & $\Sigma$ & 2 & 4 & 6 & 8 & 12 & 18 & 24 & 36 & 54 & 72 & 108 & 324 & 648 \\ \hline 
			st. trans. & 67480364637 & 67479927127 & \num{423557} & \num{12685} & 521 & 577 & 89 & 32 & 28 & 3 & 6 & 9 & 2 & 1 \\ \hline 
		\end{tabular}

		\begin{tabular}{|c|c||c|c|c|c|c|c|c|c|c|c|c|}
			\hline 
			7$\times$9 & $\Sigma$ & 2 & 4 & 6 & 8 & 12 & 18 & 24 & 36 & 54 & 108 & 216 \\ \hline 
			st. trans.      & \num{1577270689} & \num{1577056397} & \num{208362} & \num{4900} & 384 & 595 & 13 & 18 & 13 & 2 & 4 & 1 \\ \hline 
		\end{tabular}
		
		\begin{tabular}{|c|c||c|c|c|c|c|c|c|c|c|c|c|c|c|c|}
			\hline 
			8$\times$9 & $\Sigma$ & 2 & 4 & 6 & 8 & 12 & 16 & 18 & 24 & 36 & 54 & 72 & 96 & 108 & 864 \\ \hline 
			st. trans.      & \num{14940988} & \num{14917761} & \num{20900} & \num{1879} & 114 & 268 & 8 & 22 & 13 & 13 & 2 & 1 & 2 & 4 & 1 \\ \hline 
		\end{tabular}
		
		\begin{tabular}{|c|c||c|c|c|c|c|c|c|c|c|c|c|c|c|c|c|c|c|c|}
			\hline 
			9$\times$9 & $\Sigma$ & 2 & 4 & 6 & 8 & 12 & 16 & 18 & 24 & 36 & 54 & 72 & 96 & 108 & 144 & 162 & 324 & 972 & \num{7776} \\ \hline 
			st. trans.      & \num{28955} & \num{25570} & \num{2802} & 307 & 110 & 99 & 3 & 13 & 11 & 17 & 2 & 3 & 2 & 6 & 1 & 1 & 5 & 2 & 1 \\ \hline 
		\end{tabular}
		
	\caption{$2$-MOLR for $n=9$.}
	\label{9x2.}
\end{table}
\begin{table}[H]
		
		\begin{tabular}{|c|c||c|c|c|c|c|c|c|c|}
			\hline 
			2$\times$9 & $\Sigma$ & 3 & 6 & 12 & 18 & 24 & 36 & 54 & 108 \\ \hline 
			st. trans.      & 202 & 86 & 102 & 4 & 3 & 2 & 3 & 1 & 1 \\ \hline 
		\end{tabular}
		
		\begin{tabular}{|c|c||c|c|c|c|c|c|c|c|c|}
			\hline 
			3$\times$9 & $\Sigma$ & 3 & 6 & 9 & 12 & 18 & 36 & 54 & 162 & 324 \\ \hline 
			st. trans.      & \num{72836} & \num{71109} & \num{1618} & 40 & 6 & 56 & 1 & 4 & 1 & 1 \\ \hline 
		\end{tabular}
		
		\begin{tabular}{|c|c||c|c|c|c|c|c|c|}
			\hline 
			4$\times$9 & $\Sigma$ & 3 & 6 & 9 & 12 & 18 & 36 & 54 \\ \hline 
			st. trans.      & \num{1746912} & \num{1742486} & \num{4209} & 149 & 37 & 21 & 4 & 6 \\ \hline 
		\end{tabular}
		
	\caption{$3$-MOLR for $n=9$.}
	\label{9x3.}
\end{table}
\begin{table}[H]
		
		\begin{tabular}{|c|c||c|c|c|c|c|c|}
			\hline 
			2$\times$9 & $\Sigma$ & 4 &  8 & 12 & 16 & 24 & 144 \\ \hline 
            st. trans.      & \num{1017} & 881 &  121 & 5 & 2 & 7 & 1 \\ \hline 
		\end{tabular}
		
		\begin{tabular}{|c|c||c|c|c|c|c|c|c|c|c|c|}
			\hline 
			3$\times$9 & $\Sigma$ & 4 & 8 & 12 & 16 & 24 & 36 & 48 & 72 & 108 & 432 \\ \hline 
			st. trans.      & \num{356680} & \num{355023} & \num{1511} & 112 & 6 & 18 & 4 & 1 & 3 & 1 & 1 \\ \hline 
		\end{tabular}
		
		\begin{tabular}{|c|c||c|c|c|c|c|c|c|c|c|c|}
			\hline 
			4$\times$9 & $\Sigma$ & 4 & 8 & 12 & 16 & 24 & 32 & 36 & 72 & 144 & 288 \\ \hline 
			st. trans.      & \num{2640163} & \num{2635762} & \num{4131} & 95 & 147 & 1 & 17 & 2 & 5 & 2 & 1 \\ \hline 
		\end{tabular}
		
		\begin{tabular}{|c|c||c|c|c|c|c|c|c|c|c|}
			\hline 
			5$\times$9 & $\Sigma$ & 4 & 8 & 12 & 16 & 32 & 36 & 72 & 144 & 288 \\ \hline 
			st. trans.      & \num{645453} & \num{641633} & \num{3467} & 30 & 305 & 11 & 1 & 3 & 2 & 1 \\ \hline 
		\end{tabular}
		
		\begin{tabular}{|c|c||c|c|c|c|c|c|c|c|c|c|}
			\hline 
			6$\times$9 & $\Sigma$ & 4 & 8 & 12 & 16 & 24 & 36 & 48 & 72 & 216 & 432 \\ \hline 
			st. trans.      & \num{1816} & \num{1662} & 124 & 15 & 9 & 1 & 1 & 1 & 1 & 1 & 1 \\ \hline 
		\end{tabular}
		
		\begin{tabular}{|c|c||c|c|c|c|c|}
			\hline 
			7$\times$9 & $\Sigma$ & 4 & 8 & 16 & 72 & 144 \\ \hline 
			st. trans.      & 31 & 3 & 12 & 14 & 1 & 1 \\ \hline 
		\end{tabular}
		
		\begin{tabular}{|c|c||c|c|c|c|}
			\hline 
			8$\times$9 & $\Sigma$ & 8 & 64 & 288 & 576 \\ \hline 
			st. trans.      & 7 & 2 & 3 & 1 & 1 \\ \hline 
		\end{tabular}
		
		\begin{tabular}{|c|c||c|c|c|c|}
			\hline 
			9$\times$9 & $\Sigma$ & 64 & 576 & \num{2592} & \num{5184} \\ \hline 
			st. trans.      & 5 & 2 & 1 & 1 & 1 \\ \hline 
		\end{tabular}
		
	\caption{$4$-MOLR for $n=9$.}
	\label{9x4.}
\end{table}
\begin{table}[H]
		
		\begin{tabular}{|c|c||c|c|}
			\hline 
			2$\times$9 & $\Sigma$ & 5 & 10 \\ \hline 
			st. trans.      & 17 & 7 & 10 \\ \hline 
		\end{tabular}
		
	\caption{$5$-MOLR for $n=9$.}
	\label{9x5.}
\end{table}
\begin{table}[H]
		
		\begin{tabular}{|c|c||c|c|c|c|c|c|c|c|}
			\hline 
			2$\times$9 & $\Sigma$ & 6 & 12 & 24 & 36 & 54 & 72 & 108 & 216 \\ \hline 
			st. trans.      & 543 & 422 & 104 & 7 & 4 & 1 & 2 & 2 & 1 \\ \hline 
		\end{tabular}
		
		\begin{tabular}{|c|c||c|c|c|c|c|c|c|c|c|c|}
			\hline 
			3$\times$9 & $\Sigma$ & 6 & 12 & 18 & 24 & 36 & 54 & 72 & 108 & 324 & 648 \\ \hline 
			st. trans.      & \num{21620} & \num{20528} & 946 & 69 & 12 & 48 & 5 & 3 & 3 & 5 & 1 \\ \hline 
		\end{tabular}
		
		\begin{tabular}{|c|c||c|c|c|c|c|c|c|c|}
			\hline 
			4$\times$9 & $\Sigma$ & 6 & 12 & 18 & 24 & 36 & 54 & 72 & 108 \\ \hline 
			st. trans.      & 244 & 157 & 41 & 24 & 2 & 8 & 3 & 3 & 6 \\ \hline 
		\end{tabular}
		
		\begin{tabular}{|c|c||c|c|c|c|c|c|}
			\hline 
			5$\times$9 & $\Sigma$ & 6 & 12 & 18 & 36 & 54 & 108 \\ \hline 
			st. trans.      & 33 & 10 & 11 & 1 & 4 & 3 & 4 \\ \hline 
		\end{tabular}
		
		\begin{tabular}{|c|c||c|c|c|c|c|c|c|c|c|}
			\hline 
			6$\times$9 & $\Sigma$ & 6 & 12 & 24 & 36 & 72 & 108 & 216 & 324 & 648 \\ \hline 
			st. trans.      & 16 & 1 & 1 & 1 & 2 & 2 & 2 & 2 & 3 & 2 \\ \hline 
		\end{tabular}
		
		\begin{tabular}{|c|c||c|}
			\hline 
			7$\times$9 & $\Sigma$ & 12 \\ \hline 
			st. trans.      & 1 & 1 \\ \hline 
		\end{tabular}
		
		\begin{tabular}{|c|c||c|}
			\hline 
			8$\times$9 & $\Sigma$ & 48 \\ \hline 
			st. trans.      & 1 & 1 \\ \hline 
		\end{tabular}
		
		\begin{tabular}{|c|c||c|}
			\hline 
			9$\times$9 & $\Sigma$ & 432 \\ \hline 
			st. trans.      & 1 & 1 \\ \hline 
		\end{tabular}
		
	\caption{$6$-MOLR for $n=9$.}
	\label{9x6.}
\end{table}
\begin{table}[H]
		
		\begin{tabular}{|c|c||c|c|}
			\hline 
			2$\times$9 & $\Sigma$ & 7 & 14 \\ \hline 
			st. trans.      & 39 & 37 & 2 \\ \hline 
		\end{tabular}
		
		\begin{tabular}{|c|c||c|c|}
			\hline 
			3$\times$9 & $\Sigma$ & 7 & 21 \\ \hline 
			st. trans.      & \num{1532} & \num{1522} & 10 \\ \hline 
		\end{tabular}
		
		\begin{tabular}{|c|c||c|c|}
			\hline 
			4$\times$9 & $\Sigma$ & 7 & 21 \\ \hline 
			st. trans.      & 300 & 298 & 2 \\ \hline 
		\end{tabular}
		
	\caption{$7$-MOLR for $n=9$.}
	\label{9x7.}
\end{table}
\begin{table}[H]
		
		\begin{tabular}{|c|c||c|c|c|c|c|c|c|}
			\hline 
			2$\times$9 & $\Sigma$ & 8 & 16 & 32 & 48 & 56 & 96 & 864 \\ \hline 
			st. trans.      & 54 & 39 & 9 & 2 & 1 & 1 & 1 & 1 \\ \hline 
		\end{tabular}
		
		\begin{tabular}{|c|c||c|c|c|c|c|c|c|c|}
			\hline 
			3$\times$9 & $\Sigma$ & 8 & 16 & 48 & 96 & 144 & 288 & 432 & \num{2592} \\ \hline 
			st. trans.      & 48 & 22 & 17 & 3 & 2 & 1 & 1 & 1 & 1 \\ \hline 
		\end{tabular}
		
		\begin{tabular}{|c|c||c|c|c|c|c|c|c|c|c|}
			\hline 
			4$\times$9 & $\Sigma$ & 16 & 24 & 32 & 48 & 64 & 144 & 288 & 432 & 576 \\ \hline 
			st. trans.      & 27 & 11 & 2 & 2 & 3 & 4 & 1 & 1 & 1 & 2 \\ \hline 
		\end{tabular}
		
		\begin{tabular}{|c|c||c|c|c|c|c|c|c|c|}
			\hline 
			5$\times$9 & $\Sigma$ & 16 & 24 & 48 & 64 & 144 & 288 & 432 & 576 \\ \hline 
			st. trans.      & 22 & 8 & 2 & 3 & 4 & 1 & 1 & 1 & 2 \\ \hline 
		\end{tabular}
		
		\begin{tabular}{|c|c||c|c|c|c|c|c|c|c|}
			\hline 
			6$\times$9 & $\Sigma$ & 16 & 48 & 96 & 144 & 288 & 432 & 864 & \num{2592} \\ \hline 
			st. trans.      & 16 & 6 & 3 & 2 & 1 & 1 & 1 & 1 & 1 \\ \hline 
		\end{tabular}
		
		\begin{tabular}{|c|c||c|c|c|c|c|}
			\hline 
			7$\times$9 & $\Sigma$ & 16 & 48 & 96 & 288 & 864 \\ \hline 
			st. trans.      & 9 & 2 & 2 & 3 & 1 & 1 \\ \hline 
		\end{tabular}
		
		\begin{tabular}{|c|c||c|c|c|c|}
			\hline 
			8$\times$9 & $\Sigma$ & 48 & 384 & \num{1152} & \num{3456} \\ \hline 
			st. trans.      & 7 & 2 & 3 & 1 & 1 \\ \hline 
		\end{tabular}
		
		\begin{tabular}{|c|c||c|c|c|c|}
			\hline 
			9$\times$9 & $\Sigma$ & 384 & \num{3456} & \num{10368} & \num{31104} \\ \hline 
			st. trans.      & 5 & 2 & 1 & 1 & 1 \\ \hline 
		\end{tabular}
		
	\caption{$8$-MOLR for $n=9$.}
	\label{9x8.}
\end{table}
\end{landscape}

\section{The stepwise transitive $\mathbf{8}$-MOLR for $\mathbf{n=9}$}
\label{mols9}

\begin{figure}
	\begin{tabular}{|c|c|c|c|c|c|c|c|c|}
		\hline
		0 & 1 & 2 & 3 & 4 & 5 & 6 & 7 & 8 \\ \hline 
		8 & 7 & 6 & 2 & 1 & 0 & 3 & 4 & 5 \\ \hline 
		7 & 6 & 8 & 0 & 2 & 1 & 5 & 3 & 4 \\ \hline 
		6 & 8 & 7 & 1 & 0 & 2 & 4 & 5 & 3 \\ \hline 
		5 & 4 & 3 & 6 & 7 & 8 & 2 & 1 & 0 \\ \hline 
		4 & 3 & 5 & 8 & 6 & 7 & 0 & 2 & 1 \\ \hline 
		3 & 5 & 4 & 7 & 8 & 6 & 1 & 0 & 2 \\ \hline 
		2 & 0 & 1 & 4 & 5 & 3 & 7 & 8 & 6 \\ \hline 
		1 & 2 & 0 & 5 & 3 & 4 & 8 & 6 & 7 \\ \hline 
		
	\end{tabular}
	\begin{tabular}{|c|c|c|c|c|c|c|c|c|}
		\hline
		0 & 1 & 2 & 3 & 4 & 5 & 6 & 7 & 8 \\ \hline 
		7 & 6 & 8 & 0 & 2 & 1 & 5 & 3 & 4 \\ \hline 
		5 & 4 & 3 & 6 & 7 & 8 & 2 & 1 & 0 \\ \hline 
		2 & 0 & 1 & 4 & 5 & 3 & 7 & 8 & 6 \\ \hline 
		3 & 5 & 4 & 7 & 8 & 6 & 1 & 0 & 2 \\ \hline 
		1 & 2 & 0 & 5 & 3 & 4 & 8 & 6 & 7 \\ \hline 
		8 & 7 & 6 & 2 & 1 & 0 & 3 & 4 & 5 \\ \hline 
		4 & 3 & 5 & 8 & 6 & 7 & 0 & 2 & 1 \\ \hline 
		6 & 8 & 7 & 1 & 0 & 2 & 4 & 5 & 3 \\ \hline 
		
	\end{tabular}
	\begin{tabular}{|c|c|c|c|c|c|c|c|c|}
		\hline
		0 & 1 & 2 & 3 & 4 & 5 & 6 & 7 & 8 \\ \hline 
		6 & 8 & 7 & 1 & 0 & 2 & 4 & 5 & 3 \\ \hline 
		2 & 0 & 1 & 4 & 5 & 3 & 7 & 8 & 6 \\ \hline 
		3 & 5 & 4 & 7 & 8 & 6 & 1 & 0 & 2 \\ \hline 
		4 & 3 & 5 & 8 & 6 & 7 & 0 & 2 & 1 \\ \hline 
		7 & 6 & 8 & 0 & 2 & 1 & 5 & 3 & 4 \\ \hline 
		1 & 2 & 0 & 5 & 3 & 4 & 8 & 6 & 7 \\ \hline 
		8 & 7 & 6 & 2 & 1 & 0 & 3 & 4 & 5 \\ \hline 
		5 & 4 & 3 & 6 & 7 & 8 & 2 & 1 & 0 \\ \hline 
		
	\end{tabular}%
	
\caption{The stepwise transitive $8$-MOLS of size $9\times9$ 
with $|\mathrm{Aut}|=\num{10368}$, corresponding to the Galois plane.}
\end{figure}
\begin{figure}
	%
	
\caption{The stepwise transitive $8$-MOLS of size $9\times9$ with 
$|\mathrm{Aut}|=\num{31104}$, corresponding to the dual of the Hall plane.}
\end{figure}
\begin{figure}
	%
	
\caption{The stepwise transitive $8$-MOLS of size $9\times9$ with $ |\mathrm{Aut}|=384 $, corresponding to the Hall plane.}
\end{figure}
\begin{figure}
	%
	
\caption{The stepwise transitive $8$-MOLS of size $9\times9$ with $ |\mathrm{Aut}|=384$, corresponding to the dual of the Hall plane.}
\end{figure}
\begin{figure}
	%
	
\caption{The stepwise transitive $8$-MOLS of size $9\times9$ with $ |\mathrm{Aut}|=\num{3456}$, corresponding to the Hall plane.}
\end{figure}

\end{document}